\documentclass[12pt, a4paper]{article}
\usepackage{graphicx} 
\usepackage[spanish,english]{babel}
\usepackage[utf8]{inputenc}
\usepackage{amsmath} 
\usepackage{derivative}
\usepackage{tikz}
\usepackage{booktabs,caption}
\usepackage[flushleft]{threeparttable}
\usepackage{float}
\usepackage{color}
\usepackage{longtable}
\usepackage{threeparttable}
\usepackage{comment}
\usepackage{amssymb}
\usepackage{hyperref}
\usepackage{nameref}
\usepackage{csquotes} 
\usepackage[utf8]{inputenc}
\usepackage[natbibapa]{apacite}
\setcitestyle{authoryear,open={(},close={)}} 
\usepackage{lipsum} 
\usepackage[left=2.5cm, right=2.5cm, top=2.5cm, bottom=2.5cm]{geometry}
\usepackage{amsthm}
\hypersetup{
    colorlinks=true, 
    linkcolor=black, 
    citecolor=black, 
    filecolor=black, 
    urlcolor=blue, 
}
\usepackage{etoolbox}
\usepackage{array}
\newcolumntype{V}{@{}>{\hskip0pt}p{0pt}<{\vrule}}
\apptocmd{\thebibliography}{\small}{}{}
\newtheorem{theorem}{Theorem}

\newtheorem{proposition}[theorem]{Proposition}
\newtheorem{remark}[theorem]{Remark}

\title{ Evaluating the Effects of Organic vs. Conventional Farming on Aquifer Water Quality}
 \date{ }\author{ Marta Llorente$^1$ and Marta Suárez$^2$ }
\date{}

\begin{document}

\maketitle

\begin{center}
    {\small $^{1}$ Departamento de An\'{a}lisis Econ\'{o}mico: Econom\'{\i}a Cuantitativa. Universidad Aut\'{o}noma de Madrid, Campus de
Cantoblanco, 28049 Madrid, Spain.}\\
    {\small $^{2}$ Departamento de Análisis Económico Aplicado. Universidad de las Palmas de Gran Canaria, Campus Universitario de
Tafira, 35017 Gran Canaria, Spain.}\\
    {\small Emails: m.llorente@uam.es, marta.suarez@ulpgc.es}
\end{center}

\begin{abstract}	

This study analyzes water quality dynamics and aquifer recharge through irrigated agriculture, contributing to the literature on Managed Aquifer Recharge (MAR) amidst growing water scarcity. We develop two optimal control  models— a linear and a non-linear extension of  \citep{RefWorks:RefID:5-martin2013potential} -- that incorporate the impact of fertilizers on aquifer water quality, distinguishing  between organic and conventional farming practices. The linear model applies a constant rebate mechanism, whereas the non-linear model employs a concave rebate scheme.  Our results show that, depending on climate change scenarios, fertilizer-induced food price discounts,  and pollution levels, a socially optimal equilibrium in fertilizer use can be attained.  Policy implications are discussed, emphasizing the trade-off between environmental sustainability and social welfare.

\end{abstract}



\section{Introduction}

The optimal management of water resources has become a critical tool in addressing the global challenge of water scarcity. Ensuring the availability and quality of water is not only a pressing environmental issue but also a key objective of the  Sustainable Development Goals (SDGs) for the 2030 agenda.

According to the \cite{RefWorks:RefID:9-2024water}, approximately 4 billion people currently reside in water-scarce areas, and global water demand is projected to increase  by 20-30\% by 2050.  Climate change exacerbates this issue, with droughts becoming more severe. Between 2000 and 2019,  an estimated  1.43 billion people were affected by droughts. Agriculture, which accounts for 70\% of global freshwater withdrawals from rivers, lakes, and aquifers, plays a central role in this crisis. Irrigation covers 20\% of total agricultural land and contributes to 40\% of global food and fodder production. Groundwater, in particular, is a vital resource, with 2 billion people relying on it for their water supply. Effective groundwater management is therefore essential to ensure both food security and water sustainability \citep{RefWorks:RefID:9-2024water}.

In Spain, the situation mirrors global trends. Data from the  Ministry for Ecological Transition and the Demographic Challenge (MITECO) indicate that, as of 2023, 14.6\% of the country was under a state of emergency due to water scarcity, while 27.4\% was on alert \citep{RefWorks:RefID:12-miteco2023el}. Furthermore, by  May 2024, Spain's water reserves had dwindled to just  37\% of their capacity \citep{RefWorks:RefID:11-miteco2024la}. 

In this context of escalating climate uncertainty, the sustainable management of water resources —particularly concerning aquifers— has emerged as a critical area of research \citep{global_change}.  The present study addresses an optimal control problem  aimed at maximizing social welfare through the efficient allocation of water and food resources derived from aquifers. At its core, this approach emphasizes the sustainable management of these vital resources to ensure long-term societal well-being.

This study introduces two models that apply optimal control theory to groundwater management, with a focus on maximizing social welfare. These models provide a framework for evaluating policies designed to reduce fertilizer use while simultaneously enhancing social welfare within a specific aquifer system.

Building on the foundational work of \cite{RefWorks:RefID:5-martin2013potential}—which assumed zero fertilizer use, this study extends their aquifer management problem. Section \ref{sub: state} provides a detailed overview of the Martin and Stahn model, along with an in-depth exploration of the Managed Aquifer Recharge (MAR) system on which their framework is based. Expanding the scope of existing research, as outlined in Section \ref{sub: contribution}, the proposed models integrate fertilizer use and build upon the frameworks established by previous studies \citep{RefWorks:RefID:13-roseta-palma2003joint,RefWorks:RefID:14-yadav1997dynamic, RefWorks:RefID:34-brock2014optimal}. 

Fertilizer use in agriculture is known to reduce production costs, leading to lower food prices compared to organic production methods \citep{RefWorks:RefID:29-fao2024organic,agriculture11090813}.  However, the irrigation of water into aquifers through fertilizer-based crop irrigation also results in the production of polluted water. This necessitates both pumping and purification processes to render the water suitable for human consumption, incurring additional societal costs. Society derives utility from the consumption of water and food, as well as from its preferences regarding the optimal height of the aquifer and environmental quality.

To advance the analysis, this study adopts a methodological framework rooted in optimal control theory.   Section \ref{sec: Method} formalizes the class of resource management problems under investigation and establishes the necessary first- and second-order optimality conditions derived from Pontryagin’s Maximum Principle and the Arrow sufficiency theorem. These theoretical foundations provide the analytical tools to characterize solutions to the proposed models.

The subsequent sections operationalize this framework: Section \ref{sub: the problem} rigorously defines the aquifer management challenge through a bi-objective optimization lens, balancing agricultural productivity and resource sustainability. Section \ref{sub: the models} develops two distinct formulations. The first employs a linear price-discount structure to model  fertilizer-induced agricultural gains, while the second introduces a concave discount function that better captures diminishing returns and threshold effects observed in real-world systems. Both formulations explicitly incorporate hydrological constraints and socioeconomic trade-offs inherent in groundwater management.

A key assumption underpinning the models is that societal response to water pollution is influenced by nature's capacity to absorb it. High absorption capacity fosters greater societal patience, while limited capacity leads to impatience, echoing the "tragedy of the commons" articulated by \cite{tragedy}.

Theoretical solutions for the models are derived in Section \ref{sec: modelling the problem} using the Maximum Principle.  Section \ref{sub: equilibria lineal} presents a parametric and equilibria analysis of the Linear Model, identifying conditions under which fertilizer use is socially beneficial and demonstrating that total fertilizer use can be socially optimal in the absence of significant climate change impacts.  In Section \ref{sub: equilibria nolineal}, a similar analysis is conducted for the Non-Linear Model.

 Section \ref{sec: results} presents simulations for both the Linear and Non-Linear models, offering a comparative analysis to underscore key differences and insights. Section \ref{sub: lineal result} summarizes the effects of varying food rebate values in the Linear Model, while Section \ref{sub: non lin result}  examines the  Non-Linear Model, including the effects of different levels of contamination generated per unit of fertilizer.

The findings in  Section \ref{sub: discussion} reveal  a trade-off between long-term social welfare and water quality. Notably, social utility tends to be higher when fertilizers are applied at intermediate levels, as opposed to their maximum allowable quantity. To mitigate pollution, policymakers should consider measures to reduce the price gap between organic and conventional products. Such policies could facilitate a socially optimal equilibrium where water quality reaches its maximum potential. However, achieving the highest possible social utility entails accepting certain environmental costs.

This study is not without limitations, as outlined in  Section \ref{sec: Futher work}. While the developed models effectively capture global social utility, they do not account for the specific interests of key stakeholders, such as farmers, whose preferences significantly influence fertilizer usage decisions \citep{farmers, farmers_2}. Parameters used in simulations are chosen to align with \cite{RefWorks:RefID:5-martin2013potential}  ensuring comparability and enabling an evaluation of the models' performance in a generalized context. However, it is important to emphasize that these parameters are not intended to represent any specific real-world scenario, though they can be adapted for such applications. Additionally, in the Non-Linear Model, the equilibrium condition  is characterized by an implicit equation that precludes closed-form analytic solutions;  thus, we rely on linearization, Hartman’s theorem and numerical simulations to characterize the solutions and assess their stability.

The model suggests potential turnpike properties; however, the conditions required for this property to hold, as defined by \citep{marena1999neighborhood}, are not fully met. Further exploration, including alternative functional forms to refine the modeling of fertilizer benefits, is left for future research. 

In summary, Section \ref{sec: Background} introduces the problem, reviews relevant literature, and outlines this study’s contributions.  Section \ref{sec: Method} presents the methodology and theoretical framework grounded in optimal control theory.   Section \ref{sec: modelling the problem}  formulates and analytically solves the two models that incorporate fertilizer use into aquifer management. Sections \ref{sub: equilibria lineal}  and \ref{sub: equilibria nolineal} provide a parametric and equilibrium analysis of both models, identifying conditions under which fertilizer use is socially optimal. Section~\ref{sec: results} presents numerical simulations and comparative results: Section~\ref{sub: lineal result} explores the effect of varying food rebate values in the Linear Model, Section~\ref{sub: non lin result}  analyzes the Non-Linear Model under different pollution intensities, and Section~\ref{sub: discussion}  highlights the trade-off between long-term social welfare and water quality.  Section \ref{sec: Futher work} discusses the study’s limitations and outlines directions for future research. Finally,  Section \ref{sec: Conclusion}  concludes the paper.

\section{Aquifer management problem }
\label{sec: Background}

This section establishes the conceptual foundation for subsequent model development. It begins with an analysis of the optimal control framework introduced by  \citeauthor{RefWorks:RefID:5-martin2013potential}, which serves as the foundational reference for this study.

The discussion emphasizes the critical integration of water quality parameters in agricultural management models, particularly through the implementation of Managed Aquifer Recharge (MAR) systems as a mechanism for sustainable groundwater management.

This research extends existing scholarship by systematically addressing aquifer pollution challenges arising from agricultural fertilizer use. Through this analytical lens, the study offers novel insights into the water--agriculture nexus, advancing understanding of coupled human--natural systems in resource management.

\subsection{State of the art}
\label{sub: state}

 Our primary reference is \cite{RefWorks:RefID:5-martin2013potential}, which examines the Crau region in southeastern France through an optimal control framework, demonstrating how surface water irrigation maintains aquifer levels critical for both recreational ecosystems and municipal water supplies. Martin and Stahn's comparative analysis reveals superior social welfare outcomes under centralized planning (pure water policy setting) versus decentralized management by separate water and agricultural authorities (agri-environmental setting).

 While this pioneering study assumes chemical-free agriculture through a Controlled Designation of Origin framework, practical implementations frequently encounter challenges in maintaining zero-input cultivation. Furthermore, the model omits critical water quality parameters, particularly salt and nutrient accumulation from irrigation practices--a significant oversight given that salinization affects 20\% of irrigated land globally (FAO, 2021).

 Complementing this perspective, \cite{RefWorks:RefID:14-yadav1997dynamic} establishes nitrate pollution as a principal water quality constraint, while \cite{RefWorks:RefID:13-roseta-palma2003joint} develops a generalized framework for water quality management in agricultural systems.  Others like \cite{nahorski2000review}, address pollution modeling in a more general context. These works collectively underscore the necessity of dual quantity-quality optimization in groundwater management.

The methodological framework proposed in this study holds particular relevance for water-scarce regions, where Managed Aquifer Recharge (MAR) systems utilizing alternative water sources-- such as rainwater, treated wastewater, desalinated seawater-- show significant potential. Several  techno-economic analyses \citep{RefWorks:RefID:21-page2018managed, RefWorks:RefID:20-dillon2019sixty, RefWorks:RefID:17-levintal2023agricultural}  underscore the importance of MAR systems, and their economic feasibility has been examined by researchers such as \cite{RefWorks:RefID:19-heinz2011evaluating,PAVELIC2022101257}.

The models proposed in this study may also prove valuable in such contexts, as they incorporate water impurities and the associated treatment costs, which are critical factors in sustainable water management.

In the Spanish context, where  45\% of wastewater is used in agriculture, 36\% in parks and recreational areas, 10\% in industry and the remaining percentage in street cleaning and other uses \citep{RefWorks:RefID:22-jodar-abellan2019wastewater}, our model's incorporation of impurity could provide critical insights for sustainable reuse strategies. The agricultural sector is the largest recipient of wastewater, highlighting the importance of optimizing its use in this domain.

In selecting an appropriate MAR system, \cite{RefWorks:RefID:15-alam2021managed} recommend that the choice should depend on the properties, availability, and quality of the water in question. This underscores the need to incorporate water quality as a critical variable in the decision-making process.

The following section details the key aspects of the La Crau problem. For analytical consistency, we adopt the conclusions of \cite{RefWorks:RefID:5-martin2013potential}, which identify the social planner scenario as the optimal outcome.

Specifically, the social utility function employed by \cite{RefWorks:RefID:5-martin2013potential} is defined as follows:

\begin{equation}
U_0(h_t,f_t,g_t) = s(f_t) + \sigma(h_t,g_t) + e(h_t).
\label{eq:utility martin}
\end{equation}

The surplus food consumption, represented by $s(f_t)$, is defined as the difference between the willingness to pay for food consumption, represented  by the utility function $u(f_t)$, and the cost of purchasing this consumption, given by $c f_t$, where $c$ represents the unit price of food.

Similarly, the surplus groundwater consumption, denoted by $\sigma(h_t,g_t)$, is defined as the difference between the willingness to pay for the groundwater consumption and the price adjusted to the height of the aquifer, which reflects the associated pumping costs.

The final term, $e(h_t)$, captures an environmental preference. Unlike  traditional groundwater management literature  \citep{RefWorks:RefID:13-roseta-palma2003joint, JEANCHRISTOPHE2020100122},  this model assumes that consumers derive utility from the height of the aquifer, with an optimal social level denoted by $\Bar{h}$. 

Finally, the aquifer dynamics are governed by the following equation, as established by \citeauthor{RefWorks:RefID:5-martin2013potential}:

\begin{equation} \dot{h}_t = -g_t + b f_t . \label{eq: dinamic aquifer} \end{equation}

Equation \eqref{eq: dinamic aquifer}  indicates that the height of the aquifer, $h_t$, is positively influenced  by food production (through recharge mechanisms) and negatively affected  by groundwater extraction. The balance between these two factors determines the net change in aquifer height over time.

\subsection{Contribution to background literature} 
\label{sub: contribution}

Our contribution to the existing literature addresses the issue of chemicals and water quality in agriculture. Water quality is influenced by various factors, including its natural composition and contamination from human activities, such as sewage and industrial pollution. In this context, water used for agricultural irrigation--which recharges the aquifer-- can become a source of  contamination due to the application of fertilizers. Despite growing concerns regarding their environmental, health, and sustainability impacts, farmers continue to apply these chemicals due to market pressures and limited alternatives \citep{WILSON2001449}.

Furthermore, the use of fertilizers contributes to the generation of pollutants,  a process commonly assessed through nitrate concentration levels, a direct byproduct of inorganic fertilizer application on cropland \citep{RefWorks:RefID:31-bijay-singh2021fertilizers, RefWorks:RefID:30-singh2022nitrates}. Consequently, water must be treated to ensure its safety for human consumption, which introduces a social cost term into the utility function.  Equation \eqref{eq:social cost} illustrates how a decline in water quality, represented by the variable $\Psi_t$, negatively affects social welfare, following a quadratic relationship:

\begin{equation}
\rho(\Psi_t) = -\frac{1}{2} (1 - \Psi_t)^2
\label{eq:social cost}
\end{equation}
In addressing the issue of water pollution, it is important to note that the use of fertilizers—being polluting inputs—tends to lower food prices.   Therefore,  we assume that fertilizer application results in a decrease in food prices. A distinction  is drawn between organic production, which avoids fertilizer use, and conventional production, which relies on them. This choice is justified by several factors. First, organic products incur higher production costs, rendering them  more expensive than conventional ones  \citep{RefWorks:RefID:29-fao2024organic}. Second, \cite{agriculture11090813} empirically demonstrated in the  Czech Republic that conventional farming yields higher production per unit of land and lower energy consumption per unit of product compared to organic farming.

To capture this effect, we modify the functional form of the surplus from food consumption  \footnote{ $s(f_t)=u(f_t) - cf_t $} proposed  by \cite{RefWorks:RefID:5-martin2013potential}. An additional parameter, $\beta$, is introduced to represent the impact of fertilizer use on reducing food prices.  As detailed in Section~\ref{sec: modelling the problem}, two distinct modeling approaches are employed to  analyze the influence of fertilizer use on food pricing. In the first approach, the price reduction is represented as a linear function of the (dimensionless) quantity of pollutant input used at time $t$, $\gamma_t$, implying that the optimal level of fertilizer use for maximizing social welfare is either total or zero, depending on the parameter values. In the second approach, a non-linear model is considered, allowing for intermediate levels of fertilizer use and evaluating optimal social welfare under both scenarios.

These two scenarios are encapsulated in the following modification of the food surplus function:

\begin{equation}
  s(f_t, \gamma_t)  = u(f_t) - c(\gamma_t)f_t \quad \text{with } \begin{cases}
  \begin{aligned}
      &c(\gamma_t) = c - \beta \gamma_t &\quad \textbf{(Linear case)} \\
      &c(\gamma_t) = c - \beta \gamma_t^{1/2} &\quad \textbf{(Non-linear case)}
    \end{aligned}
  \end{cases}
\label{eq:fertilizer cost}
\end{equation}

For the remainder of this study, the linear representation of food price discounts resulting from fertilizer use will be referred to as the \textit{Linear Model}, while the concave representation is termed the \textit{Non-Linear Model}.

To capture the dynamics of water quality, we introduce a new equation inspired by  \cite{RefWorks:RefID:13-roseta-palma2003joint}.  This equation describes the rate of change in pollutant concentration over time:

\begin{equation}
\dot{C}_t = C(f_t, \gamma_t) - \delta C_t, \label{eq:dinamic contamination}
\end{equation} 

where  $C(f_t, \gamma_t)$ represents the concentration of contaminants generated during food production as a result of fertilizer use, $\gamma_t$ denotes  the rate of fertilizer application, and  $\delta$ corresponds to the natural degradation rate of pollutants per unit of time. Accordingly, we adopt the functional form
$$ C(f_t, \gamma_t ) = \eta f_t \gamma_t,$$
with $C_0 = 0$, where the parameter $\eta$  serves as a scaling factor  quantifying the relationship between food consumption and pollutant generation. 

Water quality is defined inversely in terms of pollution as follows:

\begin{equation}
\Psi_t = \Psi_{\text{max}} - C_t = 1 - C_t,
\label{eq:calidad agua}
\end{equation}
with the maximum water purity normalized to one. By combining the equations for pollutant concentration (\ref{eq:dinamic contamination})  and water quality (\ref{eq:calidad agua}), the evolution of water quality is expressed as:
\begin{equation}
\dot{\Psi}_t = \delta(1 - \Psi_t) - \eta f_t \gamma_t. 
\label{eq:dinamic quality water}
\end{equation}
 
This equation highlights the dual effects on water quality: while the degradation of pollutants contributes positively, the use of fertilizers has a detrimental impact.
\section{Methodology}
\label{sec: Method}

 This study extends the optimal control framework established by \cite{RefWorks:RefID:5-martin2013potential}, which demonstrates  the superiority of centralized decision-making in agri-environmental management. Our principal innovation lies in the explicit incorporation of fertilizer use as a strategic control variable, enabling analysis of its dual impact on agricultural productivity and groundwater quality.
 
The methodological approach combines analytical optimization with numerical validation. We employ Pontryagin's Maximum Principle adapted to the current-value Hamiltonian framework over an infinite time horizon (see Theorem~\ref{Maximum principle}) to derive equilibria and obtain analytical solutions for our infinite-horizon problems. Numerical simulations complement the analytical results, allowing for the characterization of dynamic trajectories. This dual methodology follows established practices in  economics focused on the efficient allocation of natural scarce resources \citep{RefWorks:RefID:33-intriligator1975applications, RefWorks:RefID:35-sethi2019optimal}.

\subsection{Optimal Control Theory over an Infinite Horizon}

The governing optimization problem follows the canonical form for stationary infinite-horizon control \citep[3.3]{RefWorks:RefID:35-sethi2019optimal}:

\begin{equation}
\left\{
\begin{aligned}
\underset{u_t \in \mathcal{U}(t)}{max} J &= \left\{ \int_{0}^{\infty} \phi(x_t,u_t)e^{-\rho t} dt \right \},\\
 &\text{subject to:}& \\
& \dot{x}_t = f(x_t,u_t), \quad x(0) = x_0,\\
\end{aligned}
\right.
\label{eq: 8 max teoria}
\end{equation}

where  $x_t \in \mathbb{R}^n$ represents the state variables, evolving according to the differential equations defined by $\dot{x}= f(x_t,u_t)$, with initial condition  $x(0)=x_0$. The control variables $u_t \in \mathcal{U}(t)\subset \mathbb{R}^m$ for $t \in (0,\infty)$, are required to be \textit{admissible} (i.e. piecewise continuous) and constrained by the set $\mathcal{U}(t)$, which  typically reflects physical or economic restrictions on the control variables at time $t$. Both state and control variables are required to be non-negative. The  functions $\phi, f  :\mathbb{R} ^n \times \mathbb{R} ^m \times \mathbb{R}  \to \mathbb{R} ^n $, are assumed to be continuously differentiable, and the  discount factor  $\rho$ lies in the range  $ [0,1]$. 

The solution framework employs the  \textit{current-value Hamiltonian}:

\begin{equation}
H(x_t,u_t,\lambda_t) = \phi(x_t,u_t) + \lambda_t f(x_t,u_t),
\label{eq: 9 (H)}
\end{equation}
where $\lambda_t$ denotes the vector of adjoint variables encoding shadow prices of resource stocks.  Under the stationary assumption, the state equation, the adjoint equation, and the Hamiltonian are explicitly independent of time $t$\footnote{Although the control and state variables are time-dependent, we use the notation $x_t$ and $u_t$ instead of the more common $x(t)$ and $u(t)$ to align with the notation used in \cite{RefWorks:RefID:5-martin2013potential}}. 

\begin{theorem}[Pontryagin’s Maximum Principle (PMP)]\label{Maximum principle} 
Given an admissible state-control trajectory  $(x^*_t, u^*_ t)$, $t \in [0 , \infty)$, that solves the optimal control 
problem  \eqref{eq: 8 max teoria}, there exists an absolutely continuous function $\lambda_t: (0, \infty)  \to \mathbb{R}^n$ such that  the following conditions are verified:

\begin{enumerate}
    \item[(i)] $\dot{x} = f(x^*_t, u^*_t, t), \quad x^*_0 = x_0, $ 
    \item[(ii)] $\dot{\lambda} = \rho \lambda_t -H_x(x^*_t, u^*_t, \lambda_t),$
    \item[(iii)] $H(x^*_t, u^*_t, \lambda_t) \geq H(x^*_t, u_t, \lambda_t )$ for each  
    $t \in [0, \infty)$ and  for all $u_t \in \mathcal {U}(t)$ .
    \item[(iv)] Transversality conditions (T.C.) \begin{equation}
\lim\limits_{t \to \infty} e^{-\rho t} \lambda_t \geq 0 \quad \text{and} \quad \lim\limits_{t \to \infty} e^{-\rho t} \lambda_t \cdot x^*_t=0.
\label{eq: tc}
\end{equation}
\end{enumerate}   
\end{theorem}
\begin{remark}
In the particular case that the solution is an interior solution,  the maximizing Hamiltonian condition (iii) is satisfied if $\frac{\partial H}{\partial u}(x_t^*,u^*_t,\lambda_t) = 0$.
\end{remark}

To ensure the optimality of a candidate ($x^*_t, u^*_t$)  that satisfies the aforementioned first-order conditions (F.O.C.), we employ \textit{Arrow's Theorem}. For this purpose, the Standard Hamiltonian is defined, where the control function  $u^*_t=u^*_t(x_t, \lambda_t)$ is implicitly and uniquely defined by:

\begin{equation}
H^0(x_t, \lambda_t) = \max_{u_t\in \mathcal{U}(t)}H(x_t,u_t, \lambda_t) =  \phi(x_t,u^*_t,t) + \lambda(t) f(x_t,u^*_t,t).
\label{eq: 9 (H_STANDAR)_second}
\end{equation}
\begin{theorem}[Sufficiency Conditions for Discounted Infinite-Horizon Problems]
Let $u^*_t$, and the corresponding $x^*_t$ and $\lambda_t$, satisfy the PMP necessary conditions  for all $t \in [0, \infty)$. Then, the pair
$(x^*_t,u^*_t)$ achieves the global maximum of  \eqref{eq: 8 max teoria} if the Standard Hamiltonian $H^0(x_t, \lambda_t)$ is concave in $x_t$ for each $t \in [0, \infty)$.
\end{theorem}

This theorem establishes that concavity of $H^0(x_t, \lambda_t)$ ensures the sufficiency of the necessary conditions for optimality, providing a guarantee that the candidate trajectory $(x^*_t,u^*_t)$ is globally optimal under the given constraints.

Finally, the optimal long-run equilibrium values for the state, control, and co-state variables, denoted as $x^{e}$, $u^{e}$ and $\lambda^{e}$\footnote{We shall use $x^{e}(1)$, $x^{e}(0)$ and $x_e$  to denote the steady state of the state variable $x_t$ in the  Linear Model with full fertilizers, the model described in \cite{RefWorks:RefID:5-martin2013potential} and the Non-Linear Model, respectively.}  must satisfy the following  conditions:

\begin{equation}
\left\{
\begin{aligned}
& f(x^{e}, u^{e})=0,\\
& \rho \lambda^{e} - H_x(x^{e},u^{e},\lambda^{e})=0,\\
& H(x^{e},u^{e},\lambda^{e}) \geq H(x^{e}, u, \lambda^{e}) \quad \forall u \in \mathcal {U}(t).
\end{aligned}
\right.
 \label{eq: teoric eq}
\end{equation}
These conditions define the steady-state behavior of the system, ensuring that the equilibrium satisfies the necessary conditions for optimality in the infinite time horizon.

\section{Modeling groundwater pollution}
\label{sec: modelling the problem}

This section applies the optimal control theory framework outlined in Section \ref{sec: Method} to model the joint dynamics of aquifer water levels and quality, and to evaluate  their impact on food and water consumption. The analysis integrates the effects of fertilizer use on these systems to derive socially optimal policies.

After defining the problem and specifying the functional forms, Section \ref{sub: the models} introduces two models, differing in their treatment of food price reductions due to fertilizer use. In Section \ref{sub: equilibria lineal}, we derive the equilibrium solutions for the Linear model and compare them to those of the  Non-Linear model in Section \ref{sub:  equilibria nolineal}. 

The Non-Linear model is introduced to address the limitations of the Linear model, as discussed in Section \ref{sub: the models}, where fertilizer use is restricted to extreme values--either complete prohibition or full application.  In contrast, the Non-Linear model allows for continuous variation in fertilizer use across its full range, thereby aiming  to maximize social welfare more effectively.

Farmer decision-making plays a crucial role in determining fertilizer use, influenced by perceived risks and incentives. For example, evidence from \citep{farmers}, indicates that farmers reduce fertilizer expenditure in response to profitability risks, while  environmental risks have little  influence on their decisions. Other factors, such as access to information, financial incentives, social norms, macroeconomic conditions, and environmental awareness, also affect these choices.\footnote{See Table 3 in \cite{farmers_2} for a review of these factors.}

Within this framework, social welfare is maximized without explicitly modeling individual agents. Instead, the environmental costs of fertilizer use-primarily related to water purification- are internalized and incorporated into the utility maximization process.

\subsection{The models}
\label{sub: the problem}

This study examines the coupled dynamics of aquifer water levels and quality, and their effects on food and water consumption. Aquifer levels depend on crop irrigation, which contributes to recharge at rate $b$, while fertilizer use  introduces contaminants at rate $\eta$. Natural pollutant degradation  occurs at rate  $\delta$. Water extracted for consumption must be treated when polluted, resulting in additional societal costs.

Consumers derive utility from the surplus of food and water consumption and exhibit environmental preferences through the desired aquifer height. When the water level deviates from the optimal benchmark,  $\bar{h}$, social utility decreases. Fertilizer use influences food prices via the cost function $c(\gamma_t)$ defined in \eqref{eq:fertilizer cost}, with price reductions increasing consumer surplus.  The parameter $d$ reflects  consumers’  willingness to pay for food.

A policy-relevant scenario  arises when $\beta<0$, corresponding to a fertilizer tax. As demonstrated in Proposition~\ref{optlinear},  fertilizer use becomes suboptimal, since the associated food price increase removes any economic incentive to apply it. This highlights the role of policy interventions in regulating fertilizer use.

The model includes three control variables: food consumption (or production) $f_t$, measured in mass  per units of time, groundwater extraction $g_t$, measured in units of volume per unit of time; and the pollutant input from fertilizer use $\gamma_t$, which is dimensionless. The state variables are the aquifer water height $h_t$, and the water quality $\Psi_t$, both dimensionless normalized between $0$ and $1$.  These together form the state vector $x(t) \in \mathbb{R}^2$ while the control vector is $u(t) \in \mathbb{R}^3$.

The dynamics of the system are  governed  by the differential equations for aquifer height and water quality, previously introduced in  \eqref{eq: dinamic aquifer}  and  \eqref{eq:dinamic quality water}, respectively, and  derived  in Section \ref{sub: contribution}. 

This model builds on the framework established in \cite{RefWorks:RefID:5-martin2013potential}, which introduces surplus functions for water consumption and environmental externalities. Our extension incorporates fertilizer use and its impact on food prices, introducing the pollutant input $\gamma_t$ directly into the surplus functions, as detailed in Section~\ref{sub: the problem}.
 
This modification enables an integrated analysis of the trade-offs between food production and environmental outcomes—specifically, the impact of fertilizer use on aquifer quality. It addresses a key limitation in the original model, which assumes chemical-free agriculture and omits pollution-related effects.

The social utility function 
\begin{equation}
U_1(h_t, \Psi_t, f_t,g_t,\gamma_t) = -\frac{1}{2}g_t^2 +g_t h_t +(d+ D(\gamma_t))f_t -\frac{1}{2} f_t^2  - \frac{1}{2} (h_t - \bar{h})^2 - \frac{1}{2} (1-\Psi_t)^2, 
\label{eq:utility final} \end{equation}
comprises three components:
\begin{itemize}
    \item \textbf{surplus from water consumption}:
    \begin{equation}
    \sigma(h_t, g_t) = -\frac{1}{2}g_t^2 + g_t h_t,
    \label{eq: exc agua}
    \end{equation}
    \item \textbf{food consumption surplus}:
    \begin{equation}
s(f_t, \gamma_t) = (d + D(\gamma_t))f_t -\frac{1}{2} f_t^2,  \text{ } \text{where}
\begin{cases}
    \begin{aligned}
    &D(\gamma_t)=\beta \gamma_t \quad \textbf{(Linear case)}  \\
    &D(\gamma_t)=\beta \gamma_t^{1/2} \quad \textbf{(Non-Linear case)},
    \end{aligned}
\end{cases} 
\label{eq: exc comida}
\end{equation}

    \item \textbf{environmental externality}:
    \begin{equation}
    e(h_t, \Psi_t) = -\frac{1}{2}(h_t - \bar{h})^2 - \frac{1}{2}(1 - \Psi_t)^2.
    \label{eq: exc envirom}
    \end{equation}
\end{itemize}

Accordingly, total social utility is given by:
\begin{equation}
U_1(h_t, \Psi_t, f_t, g_t, \gamma_t) = \sigma(h_t, g_t) + s(f_t, \gamma_t) + e(h_t, \Psi_t).
\end{equation}

The main distinction between our models and those of \cite{RefWorks:RefID:5-martin2013potential} lies in the inclusion of fertilizer-driven food price reductions through $D(\gamma_t)$ and the integration of a social cost associated with deteriorating water quality. Consequently, our social utility function  $U_1(f_t,g_t, \gamma_t, h_t, \Psi_t)$  relates to their formulation $U_0(f_t,g_t, h_t)$ (see \eqref{eq:utility martin}) as follows:

\begin{equation}
    \label{comparision}
U_1(h_t, \Psi_t, f_t,g_t, \gamma_t )= U_0(h_t, f_t,g_t)+D(\gamma_t)f_t-\frac{1}{2}(1- \Psi_t)^2.
\end{equation}

\subsection{The Problem}
\label{sub: the models}
Based on the theoretical framework established in \eqref{eq: 8 max teoria}, the problem can be formulated as the following optimal control problem:

\begin{equation}
 \left\{
\begin{aligned}
\underset{f_t, g_t, \gamma_t}{\text{max}} & \int_{0}^{\infty} \left[-\frac{1}{2}g_t^2 + g_t h_t + (d+D(\gamma_t))f_t - \frac{1}{2}f_t^2 - \frac{1}{2}(h_t - \bar{h})^2 - \frac{1}{2}(1-\Psi_t)^2\right] e^{-\rho t}  dt \\
& \text{subject to:} 
 \left\{
\begin{aligned}
 &\dot{\Psi}_t = \delta(1-\Psi_t) - \eta f_t \gamma_t, \quad  \Psi_0=1 ,\quad \Psi_t \in [0, 1],\quad \gamma_t \in [0, 1],  \\
 &\dot{h}_t = b f_t - g_t, \quad h_0 = 1, \quad h_t \in [0, 1], \\
 &\rho + \delta = 1, \quad \rho, \delta \in [0, 1], \quad \eta, b \in [0, 1],
\end{aligned}
\right.
\end{aligned}
\label{eq:optimal control}
\right .
\end{equation}
where $D(\gamma_t)$ is given by \eqref{eq: exc comida}.

The objective is to maximize the present value of social utility subject to dynamic constraints governing the evolution of water quality and aquifer height. The optimal trajectories for groundwater consumption, fertilizer use and food consumption  are determined by solving this problem under specific boundary and feasibility conditions. These trajectories reflect the decisions of a central planner seeking to achieve socially optimal outcomes.

The formulation of this problem relies on  key assumptions regarding resource sustainability and pollution dynamics. One central assumption is that the sum of the discount factor and the natural degradation rate of pollution equals one ($\rho + \delta=1$). This condition reflects the idea that future consumption is valued more when natural systems, such as aquifers, have a high capacity to purify water. 
This aligns with the notion of the "tragedy of the commons"  \citep{tragedy}, wherein over-exploitation occurs when natural resources are unable to regenerate fast enough. If the  discount rate $\rho$ approaches  its upper bound, the model predicts a decline in  groundwater consumption  toward zero, indicating low long-term sustainability. Additionally, empirical evidence from \cite{corral2002residential}  suggests that perceptions of collective behavior influence individual water consumption, with consumers being less inclined to conserve resources if they believe others are wasting them.

The model also incorporates a sustainability constraint that requires the resource growth rate $b$--representing irrigation inflows--to exceed the social discount rate $\rho$ ($b>\rho$). This condition, which builds on the findings of \cite{RefWorks:RefID:32-manzoor2014optimal}, ensures that the system can maintain long-term balance between extraction and replenishment. 

Pollution dynamics are  modeled to capture  both benefits and costs of fertilizer use. While fertilizers  improve food production, they also introduce pollutants into the aquifer. In the Linear Model, these benefits are captured through  a linear discount effect on food prices, whereas the Non-Linear Model assumes diminishing returns, consistent with studies such as \cite{RefWorks:RefID:34-brock2014optimal}. The environmental costs, primarily related to water purification, are incorporated through  convex damage functions in both cases.

For the Linear Model, the current-value Hamiltonian takes the following form:

\begin{equation}
\label{eq: current value}
\begin{aligned}
 H_1(h_t, \Psi_t, f_t,g_t,\gamma_t,\lambda_t, \mu_t)& = -\frac{1}{2}g_t^2 + g_t h_t + (d+\beta \gamma_t)f_t -\frac{1}{2} f_t^2 - \frac{1}{2} (h_t - \bar{h})^2  - \frac{1}{2} (1-\Psi_t)^2 \\
&+ \lambda_t [\delta(1-\Psi_t) -\eta f_t \gamma_t] + \mu_t [b f_t - g_t].
\end{aligned}
\end{equation}
This expression captures the dynamic interactions among state and control variables--including   aquifer height, water quality, food production, groundwater extraction, and fertilizer use--while balancing productivity gains from fertilizer use against the social costs of groundwater pollution. Since the Hamiltonian is linear in the fertilizer input $\gamma_t$, the optimal path for  $\gamma_t$ follows a bang-bang control policy. The switching rule depends on the relationship between the price-discount parameter $\beta$, the pollution factor $\eta$, and the shadow price of  water quality, $\lambda_t$:

\begin{equation}
\gamma_t^* = 
\begin{cases}
    \text{$1 \hspace{1cm} if \quad \frac{\beta}{\eta} > \lambda_t$} \\
    \text{$- \hspace{0.9cm} if \quad \frac{\beta}{\eta}= \lambda_t$} \\
    \text{$0 \hspace{1cm} if \quad \frac{\beta}{\eta} < \lambda_t$}.
\end{cases}
\label{eq: gamma lin}
\end{equation}
When  $\gamma_t^*=0$, the  optimal paths coincide with those of the baseline model  by \cite{RefWorks:RefID:5-martin2013potential}, providing a baseline for comparison against fertilizer-based scenarios.

Alternatively, the Non-Linear Model allows for intermediate levels of fertilizer application. Under certain parameter configurations, fertilizer use declines to zero when  the price reduction effect, $\beta$, diminishes  (see \eqref{eqD: demo}). The equilibrium  points \eqref{eqD: eq points} in this model are computed  numerically using the Newton-Raphson method \citep{newton}, which determines the steady-state fertilizer input $\gamma_{e}$ within the admissible range  $[0,1]$. To assess local stability, the system derived from the Maximum Principle is  linearized around the equilibrium  using Hartman’s theorem  \citep{hartmant}. 

The current-value Hamiltonian for the Non-Linear Model is given by:
\begin{equation}
\begin{aligned}
& H(h_t, \Psi_t, f_t,g_t,\gamma_t, \lambda_t, \mu_t) = -\frac{1}{2}g_t^2 +g_t h_t + (d+\beta \gamma_t^{\frac{1}{2}})f_t -\frac{1}{2} f_t^2  - \frac{1}{2} (h_t - \bar{h})^2  - \frac{1}{2} (1-\Psi_t)^2  + \\ 
 & \lambda_t [\delta(1-\Psi_t) -\eta f_t \gamma_t]  + \mu_t [b f_t -g_t].
\end{aligned}
\label{eq: current value 2}
\end{equation}

Using the respective Hamiltonians, the optimal paths are derived and validated against the transversality conditions presented in the Appendix~\ref{sec: appendix}. These trajectories offer insights into the dynamic evolution of groundwater use, food production, and pollution control. Furthermore, both models satisfy Arrow’s sufficiency conditions: their Hamiltonians are concave in the relevant state variables, guaranteeing the global optimality of the solutions under both linear and non-linear representations of fertilizer-induced food price effects.

\subsection{Analysis of Equilibria and Parameters in the Linear Model}
\label{sub: equilibria lineal}
 The linearity of the  Hamiltonian \eqref{eq: current value} in the fertilizer input yields two distinct policy regimes, determined by the shadow price of water quality,  $\lambda_t$, as compared to a threshold defined by the ratio  $\frac{\beta}{\eta}$. If the shadow price of quality water, $\lambda_t$, falls below $\frac{\beta}{\eta}$, social utility is maximized by forgoing fertilizer use. Conversely, when   $\lambda_t$ exceeds this threshold, setting $\gamma_t = 1$ becomes optimal.

We analyze the parametric constraints governing these regimes and offer a comparative analysis of equilibria under full versus zero fertilizer use, extending the framework of \citeauthor{RefWorks:RefID:5-martin2013potential}. 

We begin by establishing Proposition~\ref{optlinear}, which demonstrates that when fertilizer use leads to higher food prices ($\beta <0$), the optimal level of fertilizer application is zero ($\gamma_t^*=0$). In that case, water quality remains at its maximum level, while the optimal paths for food consumption, groundwater extraction, and aquifer height match the trajectories derived in the baseline model.

\begin{proposition}\label{optlinear}
If $\beta <0$, the optimal level of fertilizer use  in the problem \eqref{eq:optimal control} with $D(\gamma_t)=  \beta \gamma_t$ is $\gamma_t^*=0$. Additionally, the optimal path for water quality is $\psi_t^*=1$ while the optimal paths for food, groundwater consumption and aquifer height,  $( h_t^*, f_t^*, g_t^*)$ align with those  in \cite{RefWorks:RefID:5-martin2013potential}.
\end{proposition}

The proof and  corresponding expressions for $(f_t^*, g_t^*, h_t^*)$   are provided  in Appendix \ref{sub: Appendix A}. 

The conditions that lead to the optimal control $\gamma_t^*=1$  are outlined in Proposition~\ref{propgamma=1}. This analysis focuses on stable intertemporal solutions, which are critical for resource management and long-term policy planning. Stability is ensured by selecting paths associated with negative eigenvalues,  which prevent significant deviations in equilibrium due to small perturbations. Focusing the analysis on stable paths guarantees that policy implementation remains consistent over time.

When $\beta > 0$, the optimality of fertilizer use depends on the system's stability and parametric constraints. Proposition~\ref{propgamma=1} characterizes these conditions.

\begin{proposition}\label{propgamma=1}
 Let $\beta>0$ and suppose $\frac{\bar{h}}{2 - \delta}< \min \{\frac{\eta^2 + \delta(1 - b(\beta + d))}{\eta^2 + \delta(1 + b^2)}, \frac{\eta^2 + \delta - \eta(\beta + d)}{\eta b}\}$.
Then, if either $\beta > \eta$, or $\beta < \eta$ and $f^e(0)=\frac{b\bar{h}}{2 - \delta} + d \leq \frac{\beta \delta}{\eta^2}$, the stable intertemporal efficient paths for the problem \eqref{eq:optimal control} with $D(\gamma_t)=  \beta \gamma_t$  are given by \eqref{eqB: optimal path} with $\gamma_t^* = 1$.
\end{proposition}

The proof of Proposition~\ref{propgamma=1}  is detailed  in Appendix \ref{appendix B}. 

The case  $\beta > \eta$ in Proposition~\ref{propgamma=1} implies that when the reduction in food costs due to fertilizer use  exceeds the marginal environmental cost of pollution --represented in the model by the pollution coefficient $\eta$--  setting $\gamma_t^* = 1$ becomes optimal. Conversely,  if the environmental damage outweighs the cost reduction from fertilizer use ($\beta < \eta$), maximum application becomes optimal only when the combined value of water storage ($\mu^{e}(1) =  \mu^{e}(0)= \frac{\bar{h}}{\rho + 1} $) and food demand, $d$, does not exceed the benefits from fertilizer-driven food price reductions. When this condition is satisfied, the system can sustain maximum fertilizer use without compromising the resource. As evidenced in Table~\ref{tab:table 1} below, when $\beta$ is relatively low the optimal policy is to use no fertilizer, whereas for higher values  full fertilizer use is optimal. This balance highlights the role of parametric interactions in determining optimal resource management strategies.

\begin{longtable}{l r l l r r r r r r r}

\hline
Scenario   & $\beta$ & $h_e$ & $\Psi_e$ & $\gamma_e$ & $\mu_e$ & $\lambda_e$ & $g_e$ & $f_e$ & $U^e(1)$ & $U^e(0)$ \\ 
\hline
\endhead

Full  & 0.3 & 0.678 & 0.289 & 1 & 0.112 & 0.711 & 0.566 & 0.809 & 0.541 &        \\ 
Zero  & 0.3 & 0.867 &       & 0 & 0.112 &       & 0.755 & 1.078 &       & 0.589  \\ \hline
Full  & 0.4 & 0.719 & 0.237 & 1 & 0.112 & 0.763 & 0.608 & 0.868 & 0.621 &        \\ 
Zero  & 0.4 & 0.867 &       & 0 & 0.112 &       & 0.755 & 1.078 &       & 0.589  \\ \hline
Full  & 0.5 & 0.761 & 0.185 & 1 & 0.112 & 0.815 & 0.649 & 0.927 & 0.708 &        \\ 
Zero  & 0.5 & 0.867 &       & 0 & 0.112 &       & 0.755 & 1.078 &       & 0.589  \\ 
\hline
 \caption{ Comparison of equilibrium outcomes under full and zero fertilizer use for different values of $\beta$. Parameters: $\eta=0.8, \rho=0.09,   d=1, \bar{h}=0.122, b=0.7$}
 \label{tab:table 1}
\end{longtable}

A comparison of equilibria with and without fertilizer use further  highlights the nonlinear dependencies between key parameters and optimal policy outcomes.
 The Linear Model equilibria can be expressed in relation to the benchmark model proposed by \citeauthor{RefWorks:RefID:5-martin2013potential} as follows:

\begin{equation}
\begin{aligned}
&\mu^{e}(1) = \frac{\bar{h}}{\rho + 1} =  \mu^{e}(0)\\
&\lambda^{e}(1) = \frac{\eta}{(\eta^2 + \delta)}   (f^{e}(0)+ \beta) \\
& f^{e}(1) = (f^{e}(0) + \beta) \frac{\delta}{(\eta^2 + \delta)} > f^{e}(0)= d+ \frac{b\bar{h} }{1+\rho} \iff f^{e}(0) < \frac{\beta \delta}{\eta^2}\\
& \psi^{e}(1) = 1 - \lambda^{e}(1) = 1- \frac{\eta}{(\eta^2 + \delta)}   (f^{e}(0)+ \beta)  \\
& g^{e}(1)=  b f^{e}(1) > b f^{e}(0)= g^{e}(0)\iff f^{e}(0) <  \frac{\beta \delta}{\eta^2}\\ 
& h^{e}(1)= \mu^{e}(1) + g^{e}(1) > \mu^{e}(0) + g^{e}(0) =h^{e}(0)\iff f^{e}(0) < \frac{\beta \delta}{\eta^2} \\
\end{aligned}
\label{eq: comparison}
\end{equation}
Here, the subscript "$^e(0)$" refers to the equilibrium values under the benchmark model without fertilizer use ($\gamma_t^* = 0$), while "$^e(1)$" denotes the equilibrium values under full fertilizer use ($\gamma_t^* = 1$).

These relationships are summarized in the following proposition:

\begin{proposition}
For the problem \eqref{eq:optimal control} with $D(\gamma_t)=  \beta \gamma_t$, if $f^{e}(0) < \frac{\beta \delta}{\eta^2}$, the equilibrium values for food consumption, groundwater extraction, and aquifer height with full fertilizer use ($\gamma_t^* = 1$) are strictly greater than those obtained without fertilizer use ($\gamma_t^* = 0$). Conversely, if $f^{e}(0) \geq \frac{\beta \delta}{\eta^2}$, no such improvement is observed in equilibrium values.
\end{proposition}

The relationships in \eqref{eq: comparison} highlight the strong dependence of social utility on the interplay between price discounts ($\beta$), pollutant generation ($\eta$), and natural degradation ($\delta$). Utility will exceed that of the \citeauthor{RefWorks:RefID:5-martin2013potential} scenario only if price discounts and natural degradation effectively offset fertilizer-induced pollution.

\begin{figure}[H] 
    \centering 
        \centering
        \includegraphics[width=1.0\textwidth]{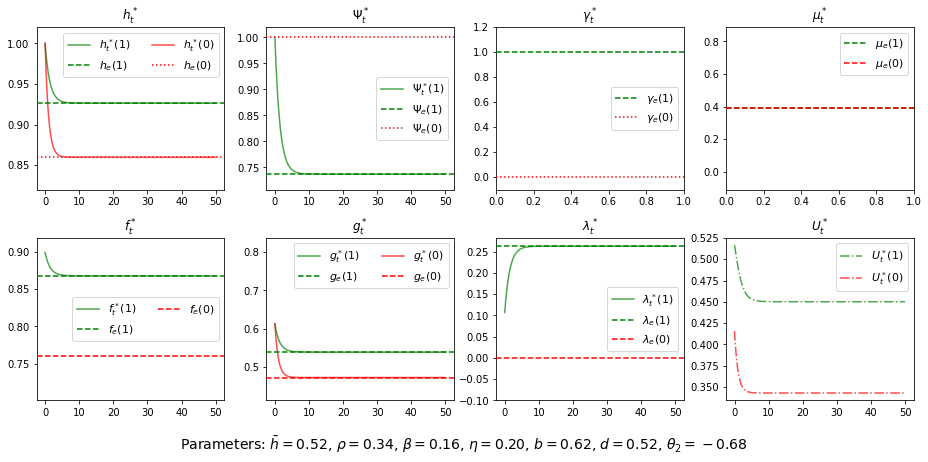} 
        \caption{In this scenario, where $f^{e}(0) \approx 0.761 < \frac{\beta \delta}{\eta^2} \approx 2.64$, full fertilizer application leads to improved equilibrium values, indicating that fertilizer-induced food price discounts, together with the system's natural degradation capacity, effectively mitigate pollutant generation.}
    \label{f0menor}
\end{figure}

\begin{figure}[H]
    \centering 
        \centering
        \includegraphics[width=1.0\textwidth]{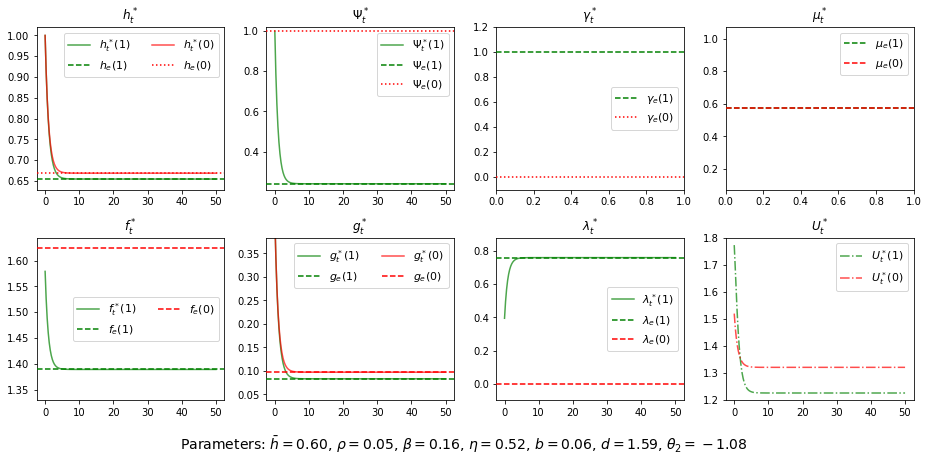} 
        \caption{In this case, we find that $f^{e}(0) \approx 1.625$, which is greater than
        $\frac{\beta \delta}{\eta^2} \approx 0.562 $. Full fertilizer application does not improve equilibrium outcomes, indicating that the benefits of reduced food prices are insufficient to offset the negative environmental impacts of pollution. 
        }
\label{f0mayor}
\end{figure}

To further illustrate these dynamics, we present two examples. Figure~\ref{f0menor} satisfies the condition $f^{e}(0) \geq \frac{\beta \delta}{\eta^2}$, while the case shown in Figure~\ref{f0mayor} does not. Both figures display the optimal trajectories for the Linear Model with and without fertilizer application. Each scenario’s corresponding equilibrium values are also shown

The shadow price $\lambda^e(1)$, associated with water quality under fertilizer use, reflects the interplay between  pollution generation, natural degradation, and food production. This balance is quantified by the ratio $\frac{\eta}{\eta^2+\delta}$. It increases with pollution ($\eta$) and decreases with stronger natural degradation ($\delta$), thus capturing the environment’s ability to offset contamination.


This trade-off has clear policy implications: robust natural degradation ($\delta$) mitigates pollution and lowers $\lambda^e(1)$. However, if pollution from fertilizers ($\eta$) far exceeds $\delta$, stringent regulatory measures may be necessary. In particular, when $\beta < 0$, fertilizer use becomes counterproductive: since $f^{e}(0) > \frac{\beta\delta}{\eta^2}$, it follows that $f^{e}(1) < f^{e}(0)$, and utility is maximized without fertilizer.

In essence, $\lambda^e(1)$ represents the marginal social cost of maintaining water quality under fertilizer use, balancing agricultural benefits against ecological damage. A high $\lambda^e(1)$  signals that pollution costs outweigh food productivity gains, justifying stricter regulations or alternative agricultural practices. Conversely, when natural degradation ($\delta$) effectively mitigates contamination, the resulting environmental damage is lower, allowing for more flexible management strategies. Ultimately, the value of $\lambda^e(1)$  reflects society’s trade-off between economic growth and environmental sustainability, guiding policy decisions on fertilizer use and water resource management.

\begin{figure}[H] \label{fig:fullbetter}
    \centering 
        \centering
        \includegraphics[width=0.8\textwidth]{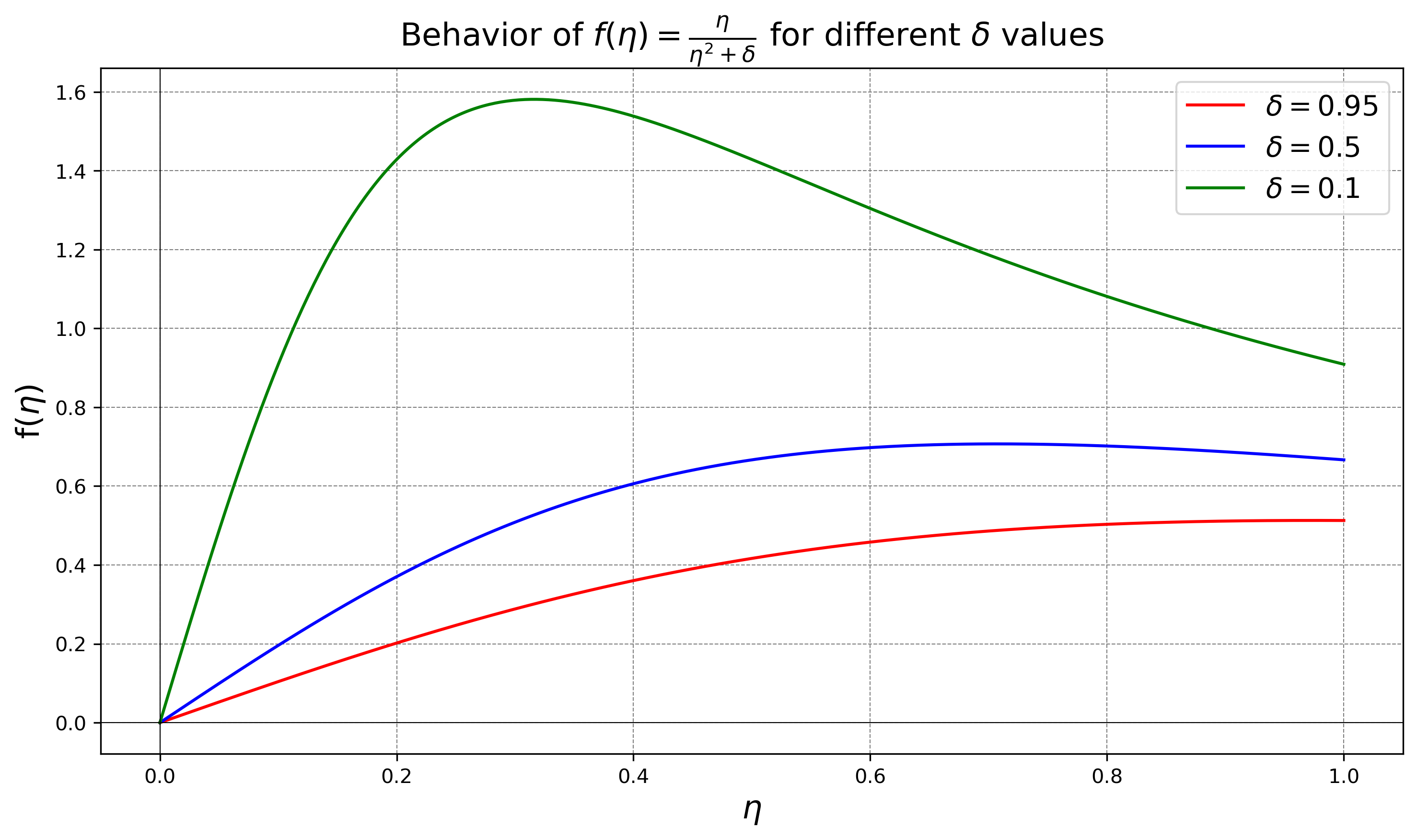} 
        \caption{If $\delta > 0.25$, then $f(\eta) < 1$ for all $\eta \in [0,1]$ while, if $0< \delta < 0.25$, $f(\eta)>1  \ \forall \eta \in [\frac{1-\sqrt{1-4\delta}}{2},\frac{1+\sqrt{1-4\delta}}{2}] $}
\label{f(eta)}
\end{figure}

\begin{figure}[H] 
    \centering 
        \centering
        \includegraphics[width=0.8\textwidth]{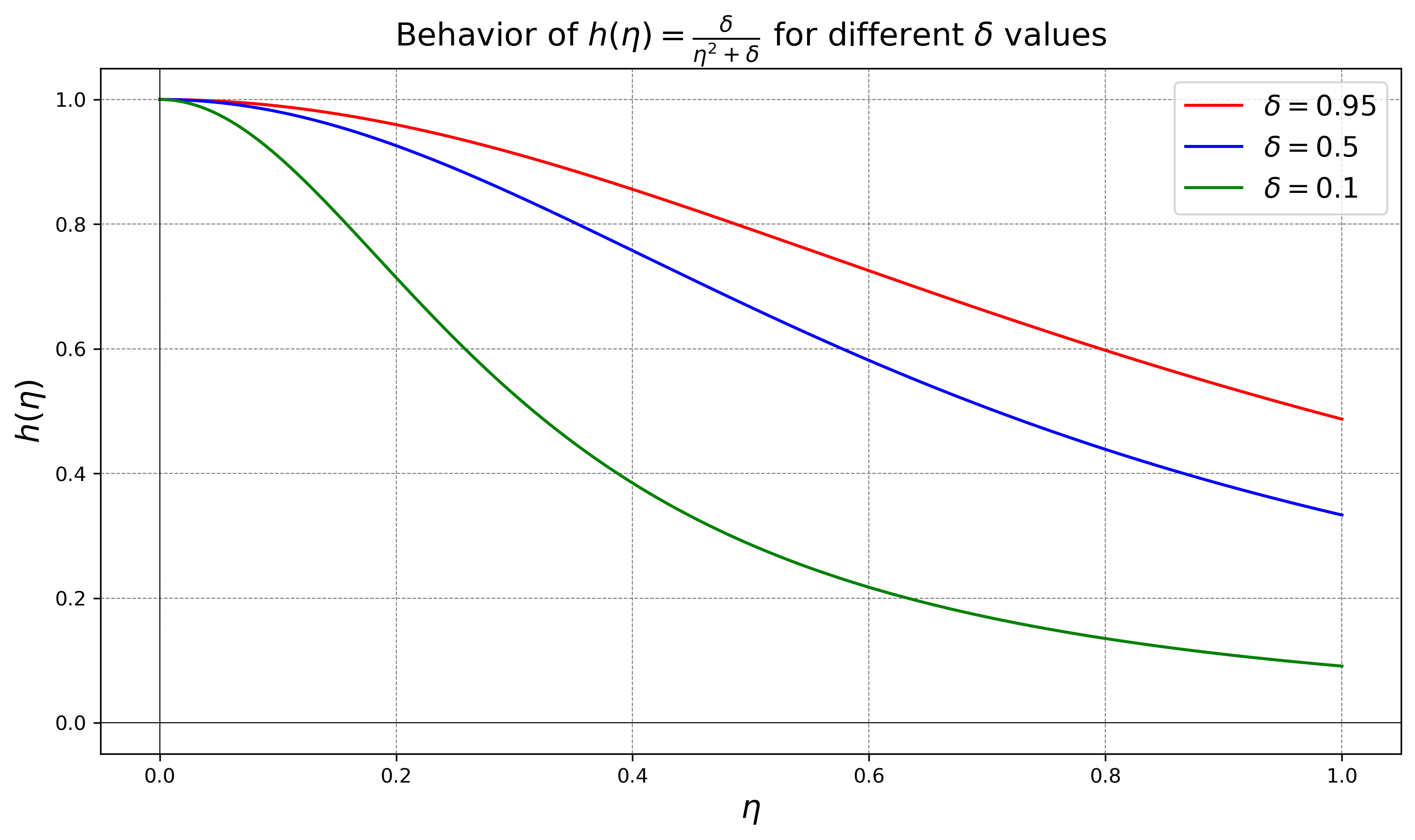} 
        \caption{} 
\label{fig: h(eta)}
\end{figure}

We now turn to the role of fertilizer-induced food price reductions in shaping system dynamics, depending on the regeneration capacity of the ecosystem and the extent of aquifer pollution. The functions $f(\eta):= \tfrac{\eta}{\eta^2 + \delta}$ and $\mathbf{h}(\eta):= \tfrac{\delta}{\delta + \eta^2}$  (illustrated in Figures~\ref{f(eta)} and \ref{fig: h(eta)}, respectively) provide further insight into the impact of pollution intensity and degradation capacity on equilibrium outcomes. Specifically, $f(\eta) $ has the property that
\begin{equation}
f(\eta) < 1 \quad \Longleftrightarrow \quad \delta > \eta - \eta^2,
\label{eq: f(eta}
\end{equation}
increasing when $\eta \in (0, \sqrt{\delta})$ and decreasing when $\eta > \sqrt{\delta}$. It achieves its maximum at $\eta = \sqrt{\delta}$, with $f(\sqrt{\delta}) = \tfrac{1}{2\sqrt{\delta}}$. If $\delta > 0.25$, then $f(\eta) < 1$ for all $\eta \in [0,1]$. Conversely, if $0 < \delta < 0.25$, there exists a range of \(\eta\) for which $f(\eta) > 1$.

Consequently, if the rate of natural degradation is high, the system can compensate for increased pollution by maintaining or improving water quality,  thus lowering $\lambda^e(1)$. A lower shadow price indicates that the system is better equipped to absorb or neutralize additional pollutant loads without compromising water quality. Conversely,  when  natural degradation is low, it becomes more difficult for the system to counteract rising pollutant levels effectively.

The function $\mathbf{h}(\eta)= \frac{\delta}{\delta + \eta^2}$ captures  how the system's capacity to degrade pollution shapes the equilibrium outcomes of  food consumption, groundwater extraction, and aquifer height. Although  $\mathbf{h}(\eta)$ does not directly represent water quality, it plays a crucial role in the system's dynamic behavior.

At $\eta =0$, $\mathbf{h}(0) =1$, suggesting that, in the absence of pollution, price discounts ($\beta$) have a stronger effect on food production. At low values of $\eta$, price discounts therefore act as an effective incentive for food production. However,  as $\eta$ increases,  $\mathbf{h}(\eta)$ decreases,  indicating that higher pollution levels reduce food production efficiency and weaken the incentive to consume water. This leads to lower water extraction and a decline in aquifer height.

When $\mathbf{h}(\eta)$  is high (e.g., due to strong natural degradation, $\delta$), the system is able to mitigate the effects of pollution more effectively, allowing higher levels of food and water consumption. Conversely, when $\mathbf{h}(\eta)$ is low --either because natural degradation is weak or pollutant input is high--pollution accumulates more rapidly.  As a result, food and water consumption must decrease to preserve long-term resource availability.

An important factor is how the discount on food prices, $\beta$, interacts with these dynamics. When $\eta$ is small, $\mathbf{h}(\eta) \approx 1$, making price discounts more effective in incentivizing increased food production. However, as pollution intensifies, this incentive weakens, leading to diminished production efficiency and reduced equilibrium resource levels.


\vspace{-2.5em} 
\begin{center}
\end{center}

This analysis illustrates how pollution input $\eta$ and the natural degradation rate $\delta$  affect the equilibrium levels of food consumption, groundwater extraction, and aquifer height. When $\delta$ is relatively high or $\eta$ is low, the equilibrium outcomes under fertilizer use improve relative to those without fertilizer, as the system offsets pollution more effectively.  Conversely, low $\delta$ and high $\eta$ reduce the beneficial effects of fertilizer-induced price discounts, potentially leading to lower equilibrium levels of food and water use.

The utility derived from full fertilizer use depends on the balance of these parameters, represented by $\frac{\beta \delta}{\eta^2}$. A higher degradation rate $\delta$ amplifies the positive effect of $\beta$ on food consumption and water use, leading to higher equilibrium utility. In contrast, when pollution generation (through large $\eta$) dominates, the environmental costs outweigh the benefits, resulting in lower social utility.

When natural purification capacity is strong,  full fertilizer use can be sustained with limited environmental risk. However, increasing pollution—often driven by climate change—may erode this resilience by overwhelming the system’s capacity to absorb contaminants. Evidence from studies on public perception of environmental issues \citep{public_perception, anderson2007exploring} suggests that public awareness can guide policy strategies to mitigate pollution and manage water resources sustainably.

Furthermore, climate change exacerbates water pollution by accelerating pollutants desorption, weakening the system’s capacity to preserve water quality. This dynamic is underscored by the United Nations:

\enquote{Water and climate change are inextricably linked. Extreme weather events are making water more scarce, more unpredictable, more polluted or all three. These impacts throughout the water cycle threaten sustainable development, biodiversity, and people’s access to water and sanitation.} \citep{ONU}.

\subsection{Analysis of Equilibria and Parameters in the Non-Linear Model}
\label{sub: equilibria nolineal}
This section analyzes the equilibrium outcomes under  a concave price discount model induced by fertilizer use, capturing  diminishing returns. We present a  key proposition characterizing  the equilibria and comparing them to both the fertilizer-free benchmark model  and the linear model with full fertilizer application.

 \begin{proposition}[Optimal equilibrium under concave fertilizer model]
\label{prop:concave_equilibrium}
Consider the concave impact of fertilizer use on food prices as described in the control problem  \eqref{eq:optimal control} with $D(\gamma_t) = \beta \gamma_t^{1/2}$. Suppose that $0<\beta \leq 2\eta$, $ 1 - \rho>\tfrac{\beta^4}{16\eta^2}$, $\tfrac{\bar{h}}{\rho+1}\le 1-\tfrac{b\beta}{2},$ along with the  following parametric constraint:
\begin{equation}
   \frac{\beta (1-\rho)}{2\eta^2} - \frac{\beta}{2} \leq b \frac{\bar{h}}{\rho+1} + d \leq \min \left\{ \frac{4 \eta (1-\rho)}{\beta^2} - \frac{\beta^2}{4\eta}, \frac{1}{b} \left(1 - \frac{\bar{h}}{\rho+1}\right) - \frac{\beta}{2} \right\}.
\label{parametric_condition}
\end{equation}

Then,  a unique equilibrium vector  $(h_e, \Psi_e,f_e, g_e, \gamma_e, \lambda_e ,\mu_e )$ exists, satisfying:
\begin{equation}
\begin{aligned}
    &\mu_e = \frac{\bar{h}}{\rho + 1}, \quad \lambda_e = \frac{\beta}{2 \eta} \gamma_e^{-1/2}, \quad \Psi_e = 1 - \lambda_e, \quad \gamma_e = \frac{\beta^2}{4 \eta^2 \lambda_e^2}, \\
    &h_e = b d + b \frac{\beta}{2} \gamma_e^{1/2} + (1+b^2) \mu_e, \quad f_e = d + \frac{\beta^2}{4 \eta \lambda_e} + b \mu_e, \quad g_e = h_e - \mu_e.
\end{aligned}
\end{equation}
where $\frac{\beta^2}{4\eta^2} \leq \gamma_e \leq 1$, and  $\Psi_e, \lambda_e \in [0, 1]$.
Moreover, the following properties hold:
\begin{enumerate}
\item If $\beta > 0 $, the equilibrium values $h_e$, $f_e$, and $g_e$ exceed those in the Linear Model without fertilizers ($\gamma_t^* = 0$).
\item If $\beta < \eta$, the equilibrium values $h_e$, $g_e$, and $f_e$  also exceed those of the Linear Model with full fertilizer use ($\gamma_t^*= 1$).
\item System stability  is ensured by the presence of two negative eigenvalues, which guarantee the  existence of stable trajectories near the equilibrium. The approximate optimal paths, derived by linearizing the system at the steady state,  are given in \eqref{eqD: optimal path}.
\item As $\beta \to 0$, fertilizer use  vanishes, and the  solution converge to those in \cite{RefWorks:RefID:5-martin2013potential}. 
\end{enumerate}

\end{proposition}
Most of the proof is given in Appendix ~ref{appendix D}; here, we focus on the equilibrium comparisons.

\begin{proof}
Let $\beta > 0$. Then,
\begin{equation}
\begin{aligned}
    &\mu_e =  \frac{\bar{h}}{\rho + 1} = \mu^e(0),  \\
    &f_e = d  + b \mu_e + \frac{\beta^2}{4 \eta \lambda_e} = f^e(0) + \frac{\beta^2}{4 \eta \lambda_e}> f^e(0) , \\ 
    &h_e = b d + (1+b^2) \mu_e+ b \frac{\beta}{2} \gamma_e^{1/2} = h^e(0)+ b \frac{\beta}{2} \gamma_e^{1/2} > h^e(0), \\ 
    &g_e = h_e - \mu_e = g^e(0)+b \frac{\beta}{2} \gamma_e^{1/2}> g^e(0).
\end{aligned}
\label{eq: comparison2}
\end{equation}
Now suppose  $\beta < \eta$. Using \eqref{eq: comparison} and \eqref{eq: comparison2} yields
\begin{equation}
\begin{aligned}
    &\mu_e = \mu^e(1) = \frac{\bar{h}}{\rho + 1},  \\
    &f_e = d  + b \mu_e + \frac{\beta^2}{4 \eta \lambda_e} =f^e(0) + \frac{\beta^2}{4 \eta \lambda_e} = \frac{\eta^2+\delta}{\delta}f^e(1)-\beta + \frac{\beta^2}{4 \eta \lambda_e}>f^e(1)+\eta -\beta>f^e(1), \\ 
    &h_e = b f_e+ \mu_e > bf^e(1)+\mu^e(1) = h^e(1)  , \\ 
    &g_e = h_e - \mu_e>h^e(1)-\mu^e(1) = g^e(1), \\ 
\end{aligned}
\end{equation}
The inequality for $f_e$ holds because $\frac{\eta^2 + \delta}{\delta} > 1$ and, since $\gamma_e \in \left[\frac{\beta^2}{4 \eta^2}, 1\right]$, it follows that implies $\lambda_e = \frac{\beta}{2 \eta} \gamma_e^{-1/2} \in \left[\frac{\beta}{2 \eta}, \frac{\beta^2}{4 \eta^2}\right]$, which ensures that $\frac{\beta^2}{4 \eta \lambda_e} > \eta$.
\end{proof}

\section{Numerical results}
\label{sec: results}
This section presents numerical simulation results, initially using the parameter set from \cite{RefWorks:RefID:5-martin2013potential}, comparing pollution models  with the benchmark model without fertilizer use.

The baseline parameters are $b=0.16, \ d=2$ and $\rho=0.05.$ It is important to note that any modifications to the values of $\eta$ and $\beta$  must be implemented consistently,  assuming all other parameters remain constant (\textit{ceteris paribus}).

The first objective of this study is to investigate the impact of the price discount parameter $\beta$ on social utility maximization and the equilibrium outcomes of the models. Subsequently, the Non-Linear Model is used to analyze the influence of both  $\beta$ and the pollutant generation parameter  $\eta$ on the optimal choice of fertilizer use. This additional analysis is particularly pertinent because, in the Linear Model with \(\gamma_t^* = 1\), fertilizer use is fixed at its maximum level, which prevents assessing the effect of varying pollution levels.

Assessing the effect of $\beta$ on social utility is critical, as this parameter may serve as a policy instrument to mitigate water pollution. For example, reducing the price discount from fertilizer use could be achieved through measures that internalize its environmental costs, such as taxes \citep{URItax199799, b6f20970-1232-3002-9edc-105d119186bc_policy}. Moreover, empirical evidence \citep{GUO2021113621_subsidie} suggests that agricultural subsidies can effectively reduce fertilizer use—and, consequently, pollution.

In Section \ref{sub: lineal result}, we summarize the effects of varying $\beta$ in the Linear Model by comparing the resulting social utility with that of the fertilizer-free benchmark. We also analyze the optimal paths for $\beta$ values that yield maximum and zero water quality, and we examine the impact of food rebates on aquifer height.

In Section \ref{sub: non lin result}, a similar analysis is performed for the Non-Linear Model, incorporating different levels of fertilizer-induced pollution. These results are compared with  the simulated food rebate scenario and the outcomes of the Linear Model. Finally, Section \ref{sub: discussion} summarizes the findings and discusses the implications of the models for policymaking.

\subsection{Simulation Results for the Linear Model}
\label{sub: lineal result}

\textbf{Simulation 1}.   Table~\ref{tab: table 2} presents the results of simulations that evaluate how varying food rebate values ($\beta$) affect the Linear Model, assuming a fixed level of pollution from fertilizer use.

\begin{longtable}{l r l l r r r r r r r}
\hline
 Scenario   &   $\beta$ &   $h_e$ &   $\psi_e$ &  $\gamma_e$ &  $\mu_e$ &   $\lambda_e$ &   $g_e$ &   $f_e$ &   $U^e(1)$ &   $U^e(0)$   \\
\hline
\endhead
 Full       &    -2.076 &   0.476 &      1.000 &            1 &     0.476 &        0.000 &  0.000 &  0.000 &     -0.000  \\  
 Zero       &    -2.076 &   0.808 &           &            0 &     0.476 &              &   0.332 &   2.076 &     & 2.163  \\  \hline 
 Full       &    -1.913 &   0.500 &      0.953 &            1 &     0.476 &         0.047 &   0.024 &   0.149 &      0.012  \\  
 Zero       &    -1.913 &   0.808 &           &            0 &     0.476 &              &   0.332 &   2.076 &      &2.163  \\  \hline 
 Full       &     0.000 &   0.780 &      0.401 &            1 &     0.476 &         0.599 &   0.303 &   1.897 &      1.967  \\  
 Zero       &     0.000 &   0.808 &           &            0 &     0.476 &              &   0.332 &   2.076 &      &2.163  \\ \hline  
 Full       &     1.000 &   0.926 &      0.113 &            1 &     0.476 &         0.887 &   0.450 &   2.810 &      4.313  \\  
 Zero       &     1.000 &   0.808 &           &            0 &     0.476 &              &   0.332 &   2.076 &     & 2.163  \\  \hline 
 Full       &     1.391 &   0.983 &     0.000 &            1 &     0.476 &         1.000 &   0.507 &   3.167 &      5.476  \\  
 Zero       &     1.391 &   0.808 &           &            0 &     0.476 &              &   0.332 &   2.076 &     & 2.163  \\  
\hline

  \caption{Equilibria for different values of $\beta$  in the Linear Model with low $ \eta$.
    Parameters: $b=0.16, d=2$, $\eta=0.3$, $\rho=0.05.$  and $\bar{h}=0.5$}
\label{tab: table 2}
\end{longtable}

\textbf{The impact of $\beta$ on  water quality and social utility}.
To achieve the highest possible water quality under fertilizer use, a sufficiently negative food price discount $\beta$ would be required. For instance, at $\beta = -2.076$, water quality reaches its upper bound ($\psi_e = 1$), but the associated social utility under full fertilizer use is slightly negative ($U^e(1) \approx -0.00028$), and thus remains well below the benchmark utility without fertilizers. Further reductions in $\beta$ lead to infeasible equilibria, as key state variables such as $\psi_t$ exceed their admissible bounds. Therefore, decreasing $\beta$ is not a viable option within the model's constraints. 

As discussed in Section~\ref{sub: equilibria lineal}, when $\beta$ is negative,  the socially optimal decision is to avoid fertilizer use entirely  ($\gamma_t^* = 0$), as this policy yields higher utility than applying fertilizers.

Table~\ref{tab: table 2} shows that increasing food price discounts ($\beta$) raises both food consumption ($f_e$) and groundwater extraction ($g_e$), resulting in higher social utility ($U^e(1)$). These improvements occur despite the deterioration in water quality ($\psi_e$), illustrating the trade-off between productivity gains and environmental costs. The results underscore how lower food prices incentivize greater resource use, especially when future environmental damages are discounted at the relatively low rate $\rho = 0.05$—a value adopted from \cite{RefWorks:RefID:5-martin2013potential} that reflects strong aquifer regeneration. This low discount rate supports a more patient societal perspective, where long-term considerations weigh heavily in maximizing social welfare.

The selection of an appropriate discount factor for environmental decisions and intergenerational equity has been widely debated for decades \citep{discount_review, discount2}. Given its critical influence on outcomes, this choice must be made carefully.  In the present framework, the discount rate directly corresponds to the aquifer's natural degradation rate ($\rho +\delta=1$).

\textbf{Impact of $\beta$ on the socially desirable aquifer height}. Under the parameters $b=0.16$, $d=2$, $\eta=0.3$ and $\rho=0.05$,  achieving the socially desirable aquifer height ($\bar{h}=0.5$) through full fertilizer use requires a specific discount value, $\beta^*=-1.9133$. As shown in Table~\ref{tab: table 2}, this choice yields an equilibrium where $h_e = \bar{h}$. Deviating from this value alters the outcome:
\begin{itemize} \item  A discount lower than $\beta^*$ leads, in equilibrium, to reduced food production, groundwater extraction, and aquifer level relative to the target $\bar{h}$.
    \item A discount above $\beta^*$ results in excess water availability  and over-incentivization of food production, as the equilibrium  aquifer level exceeds $\bar{h}$.
\end{itemize}

\begin{figure}[H] 
    \centering 
    \caption{Simulated optimal paths for the three models with $\bar{h}=0.41$, $\rho=0.04$, $\beta= 0.4$, $\eta=0.77$, $b=0.16$,  and $d=1.54$. A  high  $\eta$ significantly impacts the optimal policy, as large food discounts no longer support full fertilizer use and negatively affect water purity. Increasing  $\eta$  raises the critical $\beta$. } 
       \centering
        \includegraphics[width=0.99\textwidth]{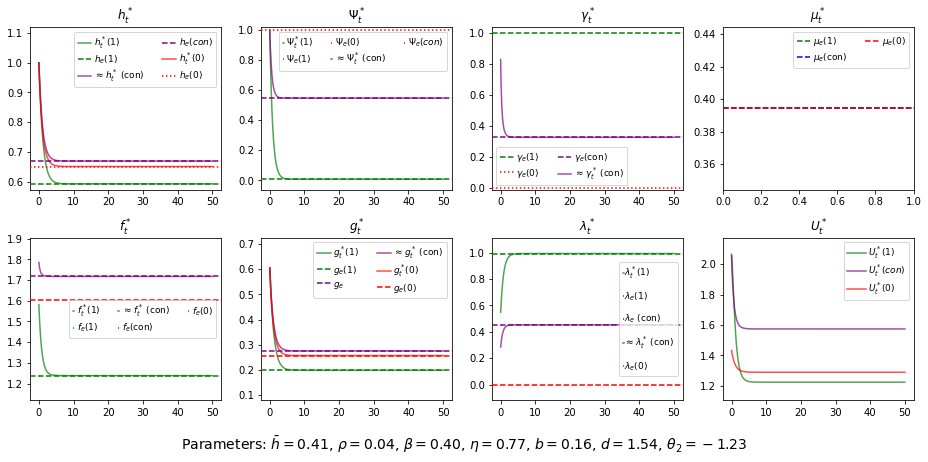} 
\label{fig:impact_eta}
\end{figure}

Furthermore, simulations in Table~\ref{tab: table 3} and Figure~\ref{fig:impact_eta} illustrate that increasing $\eta$—the rate of pollution generation—significantly alters the equilibrium outcomes. As environmental damage from fertilizer-induced contamination intensifies, a greater food price discount is required to offset its marginal social cost and validate maximum fertilizer application. This finding highlights the strong sensitivity of the optimal fertilizer policy to the parameter  $\eta$. However, this compensation has limits: for values of $\beta$ above approximately $0.435$, the model yields infeasible equilibria, as water quality ($\psi_e$) becomes negative and thus infeasible.

\begin{longtable}{l r l l r r r r r r r r}

\hline
 Scenario   &   $\beta$ &   $h_e$ &   $\psi_e$ &  $\gamma_e$ &  $\mu_e$ &   $\lambda_e$ &   $g_e$ &   $f_e$ &   $U^e(1)$ &   $U^e(0)$ \\
\hline
\endhead
 Total      &     0.000 &   0.676 &      0.214 &            1 &     0.467 &         0.786 &   0.209 &   1.044 &      0.816 &         &            \\
 Zero       &     0.000 &   0.786 &         &            0 &     0.467 &            &   0.319 &   1.593 &         &      1.279 &            \\ \hline
 Total      &     0.100 &   0.689 &      0.165 &            1 &     0.467 &         0.835 &   0.222 &   1.109 &      0.921 &         &            \\
 Zero       &     0.100 &   0.786 &         &            0 &     0.467 &            &   0.319 &   1.593 &         &      1.279 &            \\ \hline
 Total      &     0.200 &   0.702 &      0.116 &            1 &     0.467 &         0.884 &   0.235 &   1.175 &      1.033 &         &            \\
 Zero       &     0.200 &   0.786 &         &            0 &     0.467 &            &   0.319 &   1.593 &         &      1.279 &            \\ \hline
 Total      &     0.300 &   0.715 &      0.067 &            1 &     0.467 &         0.933 &   0.248 &   1.240 &      1.151 &         &            \\
 Zero       &     0.300 &   0.786 &         &            0 &     0.467 &            &   0.319 &   1.593 &         &      1.279 &            \\ \hline
 Total      &     0.400 &   0.728 &      0.017 &            1 &     0.467 &         0.983 &   0.261 &   1.306 &      1.276 &         &            \\
 Zero       &     0.400 &   0.786 &         &            0 &     0.467 &            &   0.319 &   1.593 &         &      1.279 &            \\ \hline
 Total      &     0.435 &   0.733 &      0.000 &            1 &     0.467 &         1.000 &   0.261 &    1.328 &      1.321 &         &            \\
 Zero       &     0.435 &   0.786 &         &            0 &     0.467 &            &   0.319 &   1.593 &         &      1.279 &            \\
\hline

\caption{Equilibria for different values of $\beta$ in the Linear Model with high $\eta$. Parameters: $\eta=0.7, \rho=0.07, d=1.5, \bar{h}=0.5, b=0.2$}
\label{tab: table 3}
\end{longtable}

It is worth noting that the parameter values used in Table~\ref{tab: table 3} differ slightly from those in the reference model, as they were adjusted to ensure the feasibility of the equilibria. Nonetheless, the main driver behind the observed shifts in outcomes is the increase in $\eta$.

\begin{figure}[H] 
    \centering 
    \caption{Simulated optimal paths for the three models with $\bar{h}=0.41$, \textbf{$\rho=0.70$}, $\beta= 0.4$, $\eta=0.77$, $b=0.16$,  and $d=1.54$. A lower $\delta$  significantly weakens the relative performance of fertilizer-intensive strategies. As $\delta$ increases, the environmental cost of fertilizers declines, improving the utility associated with their use. } 
       \centering
        \includegraphics[width=0.99\textwidth]{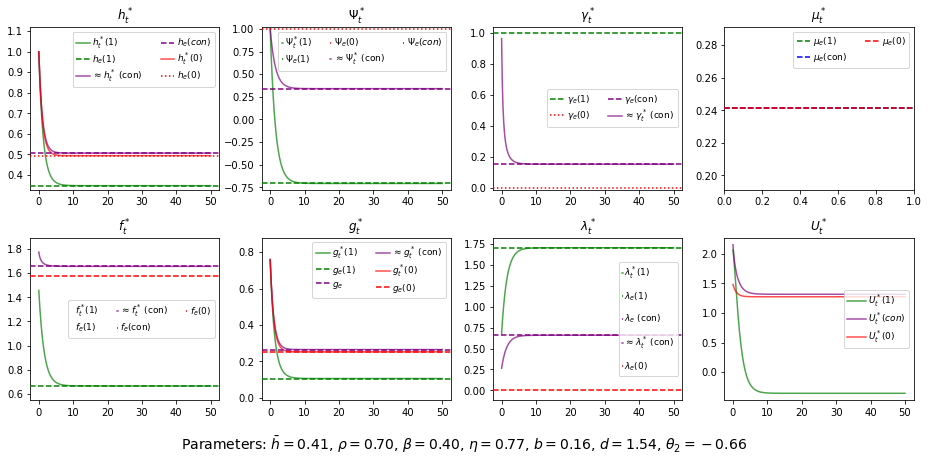} 
\label{fig:impact_delta}
\end{figure}

As illustrated in Figure~\ref{fig:impact_eta}, high levels of fertilizer-induced pollution reduce the effectiveness of food price discounts as an incentive for fertilizer use.   In such contexts, full fertilizer application ceases to be optimal. Instead, society would benefit more from either avoiding fertilizers entirely, as in the linear model, or adopting intermediate fertilizer levels that yield better utility outcomes, as captured by the non-linear model.
Moreover, additional simulations (see Figure~\ref{fig:impact_delta}) reveal that as the natural absorption capacity diminishes (i.e., when $\rho$ increases and thus $\delta = 1 - \rho$ declines), the zero-fertilizer option becomes relatively more advantageous. This underscores the role of environmental resilience in determining the desirability of fertilizer use and highlights the importance of incorporating natural degradation dynamics into policy design.

\subsection{Simulation results in the  Non-Lineal model}
\label{sub: non lin result}
\textbf{Simulation 2}. 
In the context of the Non-Linear Model, we analyze the impact of varying values of the food rebate, $\beta$, given a fixed pollution level from fertilizer use, $\eta=0.3$ (see Table~\ref{tab: table 4}). 

\begin{longtable}{l r l l r r r r r r r r}

\hline
 Scenario   &   $\beta$ &   $h_e$ &   $\psi_e$ &  $\gamma_e$ &  $\mu_e$ &   $\lambda_e$ &   $g_e$ &   $f_e$ &   $U^e(0)$ &   $U_e$ \\
\hline
\endhead
 Concave    &     0.000 &   0.808 &      0.999 &        0.001 &     0.476 &         0.001 &   0.332 &   2.076 &      &         2.163 \\
 Zero       &     0.000 &   0.808 &         &        0.000 &     0.476 &            &   0.332 &   2.076 &               2.163 &            \\ \hline
 Concave    &     0.100 &   0.813 &      0.736 &        0.397 &     0.476 &         0.264 &   0.337 &   2.108 &                 &         2.260 \\
 Zero       &     0.100 &   0.808 &         &        0.000 &     0.476 &            &   0.332 &   2.076 &               2.163 &            \\ \hline
 Concave    &     0.200 &   0.821 &      0.577 &        0.621 &     0.476 &         0.423 &   0.345 &   2.155 &                  &         2.410 \\
 Zero       &     0.200 &   0.808 &         &        0.000 &     0.476 &            &   0.332 &   2.076 &               2.163 &            \\ \hline
 Concave    &     0.300 &   0.830 &      0.441 &        0.801 &     0.476 &         0.559 &   0.354 &   2.210 &                  &         2.592 \\
 Zero       &     0.300 &   0.808 &         &        0.000 &     0.476 &            &   0.332 &   2.076 &               2.163 &            \\
\hline

\caption{Equilibria for different values of $\beta$ in the Non-Linear Model. Parameters:  $\eta=0.3, \rho=0.05, d=2, \bar{h}=0.5, b=0.16$.}
\label{tab: table 4}
\end{longtable}

\textbf{The impact of $\beta$ on fertilizers use}. Increasing the food rebate encourages greater fertilizers use in food production. This reduces water quality while raising the aquifer level due to higher food production.

It is worth noting that the equilibrium values of fertilizer use remain well below the maximum permissible level. Even when water quality is significantly degraded, the model yields moderate fertilizer application. This reflects a societal preference for more restrained agricultural practices.

\textbf{The impact of $\beta$ on quality water level and social utility}. While high food rebates negatively affect long-term water quality, lower rebates contribute positively to its preservation.

Broadly speaking, compared to the fertilizer-free equilibrium, a higher food rebate $\beta$ leads to increased food production and higher  water extraction. This is because lower food prices enhance consumer surplus, similarly to what is observed in the Linear Model.

When $\beta$ approaches zero,
the long-term utility becomes comparable to the no-fertilizer benchmark. In this scenario, aquifer height remains high, water quality approaches its maximum, and fertilizer use is neglible.

Unlike in the Linear Model, negative food discounts are not needed to discourage fertilizer use. Simply keeping rebates close to zero is enough to preserve water quality. In practical terms, this could be achieved narrowing the cost gap between organic and conventional production, for instance through technological improvements.

\textbf{Simulation 3}. We now analyze how variations in  pollutant levels, $\eta$, affect the  Non-Linear model outcomes, holding  the food rebate fixed at  $\beta=1$ (see Table~\ref{tab: table 5}).

\begin{longtable}{l r l l r r r r r r r r}
\hline
 Scenario   &   $\eta$ &   $h_e$ &   $\psi_e$ &  $\gamma_e$ &  $\mu_e$ &   $\lambda_e$ &   $g_e$ &   $f_e$ &   $U^e(0)$ &   $U^e$ \\
\hline
\endhead
 Concave    &    0.650 &   0.870 &      0.001 &        0.593 &     0.476 &         0.999 &   0.394 &   2.461 &                 &         3.487 \\
 Zero       &    0.650 &   0.808 &         &        0.000 &     0.476 &            &   0.332 &   2.076 &              2.163 &            \\
 Concave    &    0.700 &   0.867 &      0.028 &        0.540 &     0.476 &         0.972 &   0.391 &   2.444 &                         &3.420 \\
 Zero       &    0.700 &   0.808 &         &        0.000 &     0.476 &            &   0.332 &   2.076 &               2.163     &        \\
 Concave    &    0.800 &   0.862 &      0.074 &        0.456 &     0.476 &         0.926 &   0.386 &   2.414 &                 &         3.308 \\
 Zero       &    0.800 &   0.808 &         &        0.000 &     0.476 &            &   0.332 &   2.076 &               2.163    &        \\
 Concave    &    0.900 &   0.858 &      0.113 &        0.392 &     0.476 &         0.887 &   0.382 &   2.389 &                           &3.217 \\
 Zero       &    0.900 &   0.808 &         &        0.000 &     0.476 &            &   0.332 &   2.076 &               2.163     &        \\
 Concave    &    1.000 &   0.855 &      0.146 &        0.343 &     0.476 &         0.854 &   0.379 &   2.369 &                           & 3.143 \\
 Zero       &    1.000 &   0.808 &         &        0.000 &     0.476 &            &   0.332 &   2.076 &               2.163   &         \\
\hline
\caption{Equilibria for different values of $\eta$ in the Non-Linear Model. Parameters:  $\beta=1, \rho=0.05, d=2, \bar{h}=0.5, b=0.16$.}
\label{tab: table 5}
\end{longtable}

\textbf{The impact of $\eta$ on water quality and social utility}. Reducing the environmental impact of fertilizers by reducing $\eta$, encourages  fertilizer use, which in turn lowers water quality and raises the aquifer level due to increased agricultural production.

In these simulations, the Non-Linear Model consistently produces higher utility than the no-fertilizer benchmark and outperforms the cases with varying $\beta$ at fixed $\eta$.

\textbf{The aquifer level is generally above the socially desired threshold}. Simulation results consistently show that the  aquifer level exceeds the socially optimal level $\bar{h}$ across all considered cases.  This outcome highlights how even partial reductions in the pollutants generated by fertilizers--possibly achieved through technological improvements--can make fertilizer use more appealing. As a result, agricultural productivity increases, but water quality deteriorates.

\begin{figure}[H] 
    \centering
    \caption{Optimal paths in the Non-Linear Model for $\eta=0.9, \ \beta=0.85$ } 
        \includegraphics[width=1\textwidth]{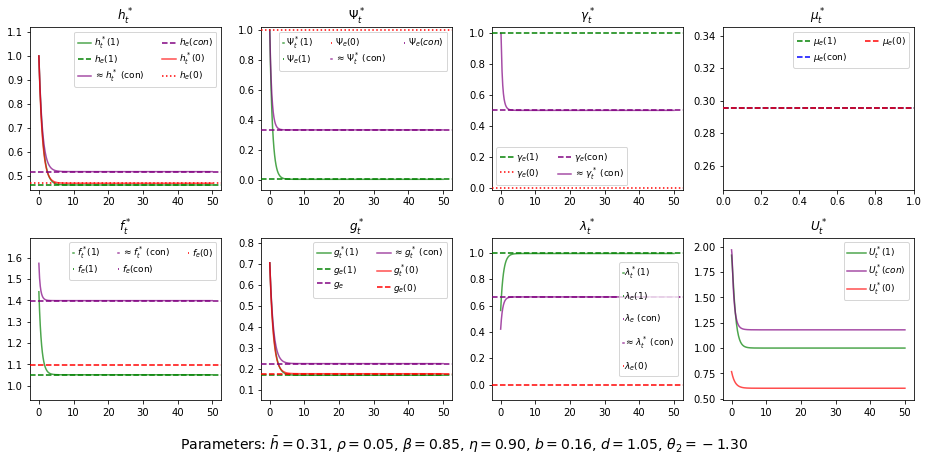} 
        \label{fig:path_etaalto_non} 
    \end{figure}

\begin{figure}[H] 
    \centering
    \caption{Optimal paths in the Non-Linear Model for $\eta=0.55,\ \beta=0.85$} 
        \includegraphics[width=1\textwidth]{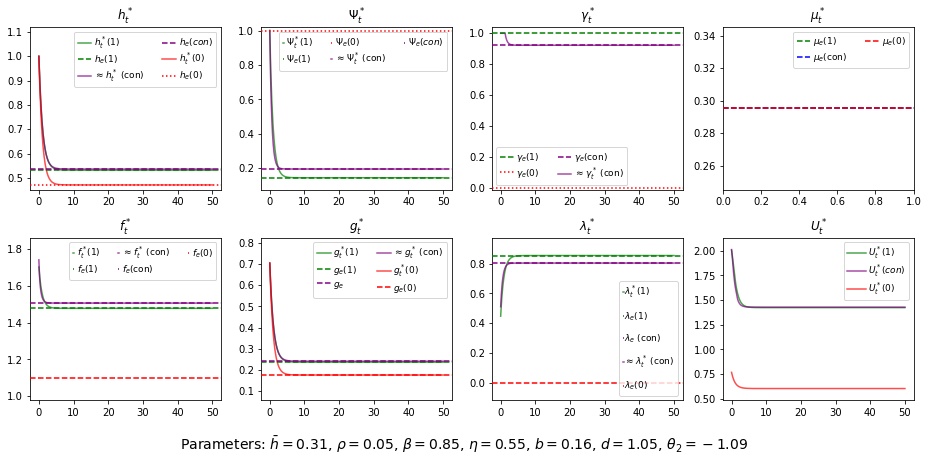} 
        \label{fig:eta_bajo_nonlineal}
\end{figure}

Figures~\ref{fig:path_etaalto_non} and~\ref{fig:eta_bajo_nonlineal} illustrate the dynamic trajectories associated with high and low $\eta$, respectively.

\vspace{-3em} 
\begin{center}
\footnotesize
   
\end{center}

Comparing  Figures \ref{fig:path_etaalto_non} and \ref{fig:eta_bajo_nonlineal} shows that social utility in both fertilizer-based models increases as the pollution coefficient $\eta$ decreases, reflecting the reduced environmental cost of fertilizer use. Furthermore, these figures illustrate that the gap between the utility trajectories in the Non-Linear and Full Fertilizer Models narrows as $\eta$ declines and, consequently, fertilizer use in the Non-Linear Model rises--driven by the lower marginal damage associated with fertilizer application. As $\gamma_t$ approaches $1$, both models converge in behavior and performance. Conversely, for higher values of $\eta$, $\gamma_e$ remains below this threshold (e.g.,  when $\eta = 0.9$, $\gamma_e \approx 0.5$), the utility  in the Non-Linear Model stays lower, though the difference progressively diminishes as $\eta$ decreases and fertilizer use intensifies.

\begin{figure}[H] 
    \centering
    \caption{Optimal paths in the Non-Linear Model for $\eta=0.55,\ \beta= 0.1$} 
        \includegraphics[width=1\textwidth]{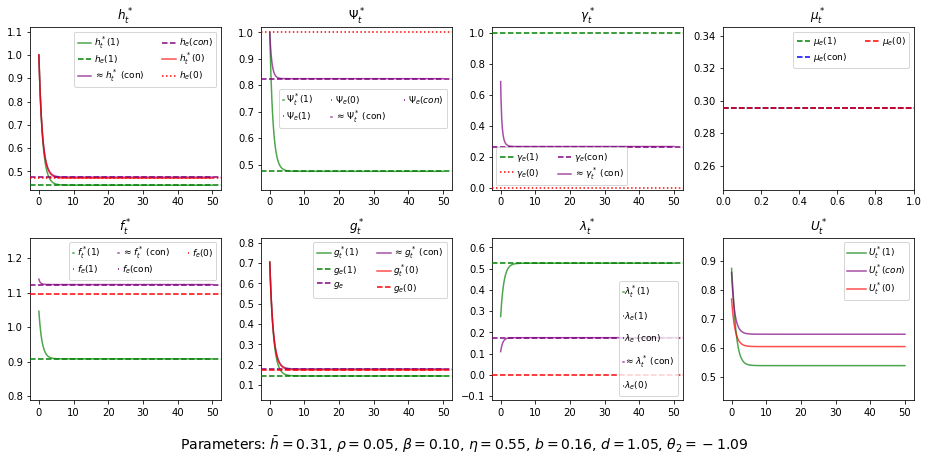} 
        \label{fig:beta_bajo_nonlineal}
\end{figure}
\vspace{-3em} 
\begin{center}
\footnotesize
   
\end{center}
An analysis of the  trajectory of $\lambda_t$ reveals how  fluctuations in fertilizer-generated pollutants or food  rebates influence the shadow price of water quality. As discussed in Section \ref{sub: equilibria lineal} this shadow price reflects the marginal cost of preserving  water quality in the presence of fertilizer use. When the pollution parameter $\eta$ decreases--as shown in Figures~ \ref{fig:path_etaalto_non} and~\ref{fig:eta_bajo_nonlineal}--the environmental cost of fertilizers is reduced, resulting in lower values of $\lambda_t$. 

In contrast, when $\beta$ increases (see Figures~\ref{fig:eta_bajo_nonlineal} and ~\ref{fig:beta_bajo_nonlineal}), the marginal cost of water quality increases, as the gains from food production no longer offset the environmental degradation. Consequently, the shadow price, $\lambda_t$, mirrors the trade-off between agricultural growth and environmental sustainability, highlighting the increasing opportunity cost of further fertilizer application.

\textbf{Further sensitivity Analysis of the Non-Linear Model.}

We present a final set of simulations examining the impact of changes in $\beta$ within the Non-Linear Model. Table~\ref{table_all} reports the values of $\beta$ for which the parametric constraints of all three models are satisfied, enabling a direct comparison. Table~\ref{table_noconcave} covers cases where the conditions of the Non-Linear Model are not fulfilled. Finally, Table~\ref{table,full_no_opt} identifies the range of $\beta$ values for which the optimal strategy involves no fertilizer use.

\begin{longtable}{lrllrrrrrrrrrV}
\hline
 Scenario   &    $\beta$ & $h_e$              & $\Psi_e$   & $\gamma_e$   & $\mu_e$             & $\lambda_e$   & $g_e$              & $f_e$                     & $U^e(1)$           &   $U^e(0)$ &   $U_e$ \\
\hline
\endhead
Concave    &     0.2 &   0.803 &      0.689 &        0.646 &     0.467 &         0.311 &   0.336 &   1.121 &       &       &       0.642 \\
 Full        &     0.2 &   0.785 &      0.545 &        1 &     0.467 &         0.455 &   0.317 &   1.058 &      0.659 &       &        \\
 No         &     0.2 &   0.779 &       &        0 &     0.467 &          &   0.312 &   1.040 &       &      0.551 &        \\ \hline
 Concave    &     0.3 &   0.820 &      0.586 &        0.819 &     0.467 &         0.414 &   0.353 &   1.176 &       &       &       0.775 \\
 Full        &     0.3 &   0.810 &      0.508 &        1 &     0.467 &         0.492 &   0.343 &   1.143 &      0.768 &      &        \\
 No         &     0.3 &   0.779 &        &        0 &     0.467 &          &   0.312 &   1.040 &       &      0.551 &        \\
\hline
\caption{Equilibria under valid parametric conditions in all three models. Parameters: $\eta=0.4, \rho=0.07, d=0.9, \bar{h}=0.5, b=0.3$.}
\label{table_all}
\end{longtable}

\begin{longtable}{lrllrrrrrrrV}
\hline
 Scenario   &    $\beta$ & $h_e$              & $\Psi_e$   & $\gamma_e$   & $\mu_e$             & $\lambda_e$   & $g_e$              & $f_e$                     & $U^e(1)$           &   $U^e(0)$    \\
\hline
\endhead
 Full        &     0.5 &   0.862 &      0.435 &            1 &     0.467 &         0.565 &   0.394 &   1.314 &      1.013 &               \\
 No         &     0.5 &   0.779 &         &            0 &     0.467 &          &   0.312 &   1.040 &       &      0.551          \\ \hline
 Full        &     0.6 &   0.887 &      0.398 &            1 &     0.467 &         0.602 &   0.420 &   1.399 &      1.148 &                \\
 No         &     0.6 &   0.779 &        &            0 &     0.467 &         0.000 &   0.312 &   1.040 &       &      0.551         \\ \hline
 Full        &     0.7 &   0.913 &      0.361 &            1 &     0.467 &         0.639 &   0.445 &   1.485 &      1.292 &               \\
 No         &     0.7 &   0.779 &         &            0 &     0.467 &         0.000 &   0.312 &   1.040 &       &      0.551         \\ \hline
 Full        &     0.8 &   0.938 &      0.325 &            1 &     0.467 &         0.675 &   0.471 &   1.570 &      1.444 &                \\
 No         &     0.8 &   0.779 &        &            0 &     0.467 &         0.000 &   0.312 &   1.040 &      &      0.551          \\ \hline
 Full        &     0.9 &   0.964 &      0.288 &            1 &     0.467 &         0.712 &   0.497 &   1.655 &      1.604 &               \\
 No         &     0.9 &   0.779 &        &            0 &     0.467 &         0.000 &   0.312 &   1.040 &      &      0.551          \\ \hline
 Full        &     1 &   0.990 &      0.251 &            1 &     0.467 &         0.749 &   0.522 &   1.741 &      1.773 &              \\
 No         &     1 &   0.779 &       &            0 &     0.467 &         0.000 &   0.312 &   1.040 &       &      0.551         \\
\hline
\caption{Equilibria when the Non-Linear Model’s parametric restrictions are not satisfied. Parameters: $\eta=0.4, \rho=0.07, d=0.9, \bar{h}=0.5, b=0.3$. }
\label{table_noconcave}
\end{longtable}

\begin{longtable}{lrllrrrrrrV}
\hline
 Scenario   &    $\beta$ & $h_e$                & $\gamma_e$   & $\mu_e$               & $g_e$              & $f_e$                     & $U^e(1)$           &   $U^e(0)$ \\
\hline
\endhead

 No         &     0.050 &   0.779         &            0 &     0.467 &              0.312 &   1.040 &      0.510 &      0.551          \\
 No         &     0.000 &   0.779         &            0 &     0.467 &              0.312 &   1.040 &      0.465 &      0.551         \\
 No         &     0.160 &   0.779         &            0 &     0.467 &              0.312 &   1.040 &      0.617 &      0.551         \\
 No         &     0.400 &   0.779         &            0 &     0.467 &               0.312 &   1.040 &      0.887 &      0.551        \\
\hline
\caption{Equilibria where full fertilizer use is suboptimal. Parameters: $\eta=0.4, \rho=0.07, d=0.9, \bar{h}=0.5, b=0.3$. }
\label{table,full_no_opt}
\end{longtable}

\subsection{Results discussion}
\label{sub: discussion}

The analysis conducted in this paper demonstrates the usefulness of deterministic dynamic models for evaluating the long-term impact of fertilizer use on aquifer systems. By incorporating key environmental and economic parameters, the models enable a comprehensive assessment of optimal societal decisions concerning food production, water consumption, and fertilizer application, along with their effects on aquifer water quality and depth.

The two proposed models—linear and non-linear—enable a structured comparison between three different scenarios: (i) no fertilizer use, following the benchmark model of \cite{RefWorks:RefID:5-martin2013potential}; (ii) full fertilizer application; and (iii) intermediate fertilizer use levels resulting from concave price discounting. These models are flexible enough to identify under which combinations of parameters each scenario emerges as optimal in terms of social utility.

Moreover, the framework allows for the exploration of how parametric variations influence model outcomes, offering insights into plausible real-world dynamics in aquifers with similar characteristics. These simulations facilitate the design and evaluation of targeted environmental and agricultural policies, such as reducing price incentives for fertilizer use or mitigating the pollutants generated by agricultural practices. The outcomes shed light on how policy choices shape optimal societal behavior and, ultimately, the trade-off between environmental preservation and economic efficiency.

With respect to the conditions under which fertilizer use is not socially desirable, the Linear Model shows that negative  price rebates (or even moderate positive when the environmental impact of fertilizers is high) are necessary to discourage fertilizer application. The Linear Model also underscores the role of environmental resilience in shaping long-term outcomes. As the natural absorption capacity $\delta$ decreases, the relative difference in utility between the fertilizer-based and the no-fertilizer models grows, with the latter becoming comparatively more advantageous. This highlights the importance of environmental degradation rates in evaluating policy alternatives.

In the Non-Linear Model, numerical simulations evidence that  increasing  price discounts and  decreasing pollution per unit of fertilizer contribute to higher fertilizer use. However, even when water quality deteriorates substantially, the optimal level fertilizer remain moderate--well bellow the maximum threshold. This reflects a built-in precautionary behavior in the model structure, offering a more realistic depiction of decision-making under environmental constraints. 

Importantly, the Non-Linear Model outperforms the Linear Model with full fertilizer use in intermediate scenarios, where neither the pollution load nor the price discount is extreme. In these contexts, moderate fertilizer use maximizes social welfare more effectively than the two polar alternatives.

If the primary goal is to preserve water quality, minimizing food price discounts associated with fertilizer use is more effective than attempting to reduce the pollutants themselves. The latter strategy tends to incentivize fertilizer use, thus increasing environmental pressure.

The models suggest that society exhibits a higher elasticity in pollutant-fertilizer consumption, meaning that increases or decreases in pollutants lead to proportional changes in fertilizer use. Conversely, the impact of food price rebates is more dependent on pollution levels, requiring larger discounts to incentivize fertilizer use when pollution is high.

Finally, from a policy standpoint, raising awareness of the environmental externalities of fertilizers, and promoting innovations that reduce their harmful effects, appears to be a more efficient strategy than providing price incentives for conventional agricultural production.

In summary, our results indicate that intermediate fertilizer use consistently leads to higher social welfare than  zero or full fertilizer application. However, it is important to emphasize that these results are grounded in numerical simulation and contingent on the specific assumptions and parameter values used.

\section{Limitations and Further work}
\label{sec: Futher work}
This study has several limitations inherent to its theoretical and numerical nature. The models developed focus exclusively on fertilizer use in agriculture, providing a framework for analyzing the long-term impacts of such practices on aquifer systems. While the models yield valuable insights into policy evaluation and water quality dynamics, the conclusions are conditional on the model assumptions and parameter values.
 
For instance, the representation of social utility in the model adopts a generic formulation that does not explicitly account for stakeholder-specific preferences, such as those of farmers, whose objectives may differ from those of society as a whole~\citep{farmers, farmers_2}.

To facilitate a direct comparison with the benchmark model by \citeauthor{RefWorks:RefID:5-martin2013potential}~\citep{RefWorks:RefID:5-martin2013potential}, the parameter values used here were deliberately chosen to align with or remain close to those in the reference framework. While this enhances comparability, future work should incorporate empirical data to validate the results and improve the model’s applicability to real-world settings.

For example, parameter calibration based on observed data from a specific geographical area could provide more realistic insights. In addition, the sensitivity of model dynamics to variations in natural degradation capacity deserves further exploration, particularly in light of the assumptions made regarding the discount rate.

Finally, it is also important to note that the current models attribute aquifer contamination exclusively to fertilizer use. In practice, groundwater pollution is shaped by a broader set of factors, including saline intrusion, wastewater discharge, and soil aridity~\citep{global_change}. Expanding the model to incorporate these elements, as well as Managed Aquifer Recharge (MAR) strategies, could offer a more comprehensive understanding of water quality dynamics.

\section{Conclusions}
\label{sec: Conclusion}

This paper has addressed the management of groundwater resources in the context of fertilizer-dependent agriculture, with particular attention to its impact on aquifer water quality and depth. To this end, two dynamic models based on optimal control theory were developed, building on the benchmark framework proposed by \citeauthor{RefWorks:RefID:5-martin2013potential}~\citep{RefWorks:RefID:5-martin2013potential}.

The study contributes to the literature on groundwater pollution by modeling fertilizer use as both a source of environmental degradation and a mechanism for lowering food prices through reduced production costs. While fertilizers may offer economic advantages over organic alternatives, they also entail significant environmental costs due to water contamination, which necessitates treatment and management.

The two models differ in how they represent the cost-reduction effect of fertilizers. The Linear Model assumes an all-or-nothing approach to fertilizer use, while the Non-Linear Model introduces a concave discount structure that allows for interior optima. This distinction enables a richer analysis of policy scenarios and equilibrium outcomes.
 
The analytical and numerical resolution of the models, based on the Maximum Principle and complemented by Hartman's linearization, has allowed for a precise characterization of equilibrium behavior under different policy and environmental scenarios.

A series of simulations based on the parameter set from \citeauthor{RefWorks:RefID:5-martin2013potential} was carried out to evaluate how changes in price discounts and pollution levels affect system dynamics and social welfare. Simulation 1 explored the effects of fertilizer-induced food price reductions in the Linear Models. Simulations 2 and 3 focused on the Non-Linear Model, analyzing the optimal fertilizer levels under varying discount rates and the effects of changes in pollutant intensity per unit of fertilizer applied.

The analysis of the Linear Model revealed the conditions under which full or zero fertilizer use emerges as optimal. These outcomes depend critically on the pollution level, the natural degradation capacity, and the strength of price incentives. When natural absorption is high and pollution is low, full fertilizer use can be socially desirable. Conversely, steep price discounts may promote fertilizer use even when it is environmentally detrimental.

The results of the Non-Linear Model  indicate  that intermediate levels of fertilizer use can yield  higher social utility than either  of the two extremes. Even under substantial deterioration of water quality, optimal fertilizer application remains below the maximum threshold, suggesting a precautionary behavior embedded in the model.

To address the challenges posed by climate change, it is advisable to implement policies aimed at reducing the price differential between organic and conventional agricultural practices. Narrowing this gap would help curb fertilizer use by making sustainable alternatives more economically viable. One promising avenue is investment in research and development to lower the production costs associated with non-polluting methods. This is particularly relevant given that organic products tend to be more expensive due to their inherently higher production costs \citep{WILSON2001449, RefWorks:RefID:29-fao2024organic, agriculture11090813}.

Alternatively, targeting the pollutant intensity of fertilizer use may prove more effective than modifying price incentives. The models suggest that fertilizer use is more responsive to changes in pollution levels than to changes in food price discounts. Thus, mitigating environmental damage at the source could yield greater improvements in welfare while exerting less pressure on water resources.

Ultimately, the proposed models generalize the benchmark case by \citeauthor{RefWorks:RefID:5-martin2013potential} and replicate its results when no fertilizers are applied. They offer a versatile tool for evaluating environmental and agricultural policies aimed at balancing sustainability and economic efficiency in aquifer systems.

While the models show promise for policy analysis in aquifers with similar characteristics, it is important to acknowledge their limitations. The findings are contingent on the theoretical structure and parameter assumptions, and may not generalize without empirical validation. Future research should expand the empirical basis of the models, conduct more comprehensive parametric analyses, and explore the existence of turnpike properties to strengthen the robustness of the conclusions.
\section*{Acknowledgements}
Marta Llorente acknowledge ﬁnancial support from UCM-Santander project PR3/23-30830.

\clearpage \bibliography{export}


\section{Appendix}
\label{sec: appendix}

\setcounter{equation}{0}
\setcounter{subsection}{0}
\renewcommand{\thesubsection}{\Alph{subsection}} 
\renewcommand{\theequation}{\thesubsection.\arabic{equation}} 
\renewcommand{\thefigure}{\thesubsection.\arabic{figure}} 
\renewcommand{\thetable}{\thesubsection.\arabic{table}} 

\subsection{Proof of Proposition~\ref{optlinear}}
\label{sub: Appendix A}
\begin{proof}
To verify the first-order conditions specified  by the Pontryagin’s Maximum Principle (see Theorem~\ref{Maximum principle}), we define the current-value Hamiltonian for the linear problem \eqref{eq:optimal control} with $D(\gamma_t)= \beta \gamma_t$ as:
\begin{equation*}
  H_1(h_t, \Psi_t,f_t,g_t,\gamma_t, \lambda_t, \mu_t)=  U_1(h_t, \Psi_t,f_t,g_t,\gamma_t) + \mu_t( b f_t - g_t)+\lambda_t [\delta(1 -\Psi_t) - \eta f_t\gamma_t].
\end{equation*} 
Using \eqref{comparision}, this expression can be rewritten as:
\begin{eqnarray}\label{Hamiltoniansrel}
 H_1(h_t, \Psi_t,f_t,g_t,\gamma_t, \lambda_t, \mu_t)&=&U_0(h_t,f_t,g_t) +\beta \gamma_tf_t-\frac{1}{2}(1-\Psi_t)^2 + \mu_t( b f_t - g_t) \nonumber\\
&+&\lambda_t [\delta(1 -\Psi_t) - \eta f_t\gamma_t]=\\
&=& H_0(h_t,f_t,g_t,\mu_t) +\beta \gamma_ tf_t-\frac{1}{2}(1-\Psi_t)^2 +\lambda_t [\delta(1 -\Psi_t) - \eta f_t\gamma_t],\nonumber
\end{eqnarray}
where 
\begin{equation*}
  H_0(h_t,f_t,g_t,\mu_t)=  U_0(h_t,f_t,g_t) + \mu_t( b f_t - g_t)
\end{equation*}
is the current value Hamiltonian for  the benchmark problem:

\begin{equation}
 \left\{
\begin{aligned}
\underset{f_t, g_t}{\text{max}} & \int_{0}^{\infty} U_0(h_t,f_t,g_t) e^{-\rho t}  dt \quad   \rho \in [0, 1]\\
& \text{s.t.} \quad  \dot{h}_t = b f_t - g_t, \quad h_0 = 1, \quad h_t \in (0, 1). \\
\end{aligned}
\label{Martinproblem}
\right.
\end{equation}
The optimal intertemporal paths of problem \eqref{Martinproblem}, as derived in \cite{RefWorks:RefID:5-martin2013potential}, are given by:
\begin{equation}
\begin{aligned}
f_t^* &= d + \frac{b\Bar{h}}{1+ \rho}, \quad  g_t^*= g^{e}(0) +(1-h^{e}(0))  e^{-t} ,\\
h_t^*&=h^{e}(0)+(1-h^{e}(0))e^{- t}, \quad \mu^*_t=\frac{\bar{h}}{\rho +1},\\
\end{aligned}
\label{Martinsol}
\end{equation}
where: $$g^{e}(0)=bd+\frac{\bar{h}b^2}{\rho +1} \ \text{ and  } h^{e}(0)=bd+\frac{\bar{h}(b^2+1)}{\rho +1},$$

From Theorem~\ref{Maximum principle}-(iii), it follows that:
\begin{eqnarray}\label{Maxpiii}
H_0(h^*_t,f_t^*,g_t^*, \mu_t^*)\ge  H_0(h^*_t, f_t,g_t, \mu_t^*) 
\end{eqnarray}
for any admissible control pair $(f_t,g_t)$.
Now, by setting $\gamma_t^* = 0, \ \Psi_t^*= 1 , \ \lambda_t^*=0$, we obtain
\begin{eqnarray}\label{Hamiltoniansrel2}
H_1(h^*_t, \Psi^*_t,f_t^*,g_t^*,\gamma^* _t,\lambda_t^*, \mu_t^*)
= H_0( h^*_t, f_t^*,g_t^*,\mu_t^*).
\end{eqnarray}
 
Suppose that $\beta<0$, combining \eqref{Hamiltoniansrel2}, \eqref{Maxpiii} and \eqref{Hamiltoniansrel} yields:
\begin{eqnarray*}
H_1( h^*_t, \Psi^*_t, f_t^*,g_t^*,\gamma^* _t,\lambda_t^*, \mu_t^*)&=&H_0(h^*_t,f_t^*,g_t^*, \mu_t^*)\ge  H_0(h^*_t, f_t,g_t, \mu_t^*) 
= \\
&=& H_1( h^*_t, \Psi^*_t, f_t,g_t,\gamma_t,\lambda_t^*, \mu_t^*)-\beta \gamma_t f_t \ge H_1(h^*_t, \Psi^*_t, f_t,g_t,\gamma _t, \lambda_t^*, \mu_t^*)  ,
\end{eqnarray*}
which proves that the controls $(h^*_t, \Psi^*_t)$ satisfy condition (iii) of Theorem~\ref{Maximum principle} for the adjoint variables $\lambda_t^*$ and  $\mu_t^*$.

Condition (i) holds trivially. To verify condition (ii), we use the optimality of the paths $(h^*_t, f_t^*,g_t^*, \mu_t^*)$ and the relation in  \eqref{Hamiltoniansrel}, to obtain 
$$\dot{\mu}^*_t= \rho \mu_t^* - \frac{\partial H_0(h^*_t,f_t^*,g_t^*,  \mu_t^*)}{\partial h_t}= 
\rho \mu_t^* - \frac{\partial H_1( h^*_t, \Psi^*_t,f_t^*,g_t^*,\gamma^*_t, \lambda_t^*, \mu_t^*)}{\partial h_t}.$$
The equation for $\lambda_t^*$ holds naturally as:
$$\dot{\lambda}^*_t= \rho \lambda_t^* - \frac{\partial H_1( h^*_t, \Psi^*_t,f_t^*,g_t^*,\gamma^*_t, \lambda_t^*, \mu_t^*)}{\partial \Psi_t}=0,$$
given that $\lambda_t^*=0 $ and $\Psi^*_t=1.$

The transversality conditions 
hold trivially since $\mu_t$ is constant, $\rho <1$ and $\lambda_t=0$.

Finally,  the second-order conditions are fulfilled  because the Standart Hamiltonian 
\begin{equation}
\label{eqA: current value second order}
\begin{aligned}
&  \hat{H}( h_t, \Psi_t, \lambda_t, \mu_t) = e^{-t\rho} (-\frac{1}{2}{g_t^*}^2 + g_t^* h_t + df_t^* -\frac{1}{2} {f_t^*}^2 - \frac{1}{2} (h_t - \bar{h})^2  - \frac{1}{2} (1-\Psi_t)^2) \\
&+ \lambda_t \delta(1-\Psi_t)+ \mu_t [b f_t^* - g_t^*]
\end{aligned}
\end{equation}
is concave in $( h_t, \Psi_t).$
   
\end{proof}

\subsection{Proof of Proposition~\ref{propgamma=1}}
\label{appendix B}
\setcounter{equation}{0}
\begin{proof}
Let $\beta >0$ and suppose that
\begin{equation}
 \label{linealparamrestric}
 \mu^{e}=\frac{\bar{h}}{2- \delta}< \min\{m_1=:\frac{\eta^2+\delta -\delta b(\beta +d))}{\eta ^2+\delta(1+b^2)}, m_2=:\frac{\eta ^2+\delta-\eta(\beta+d)}{\eta b}\}.
\end{equation}
Define the current-value Hamiltonian of problem \eqref{eq:optimal control} with $D(\gamma_t)= \beta \gamma_t$ as:
\begin{eqnarray*}
  H_1(h_t, \Psi_t,f_t,g_t,\gamma_t,  \lambda_t, \mu_t)&=&  -\frac{1}{2}{g_t}^2 + g_t h_t + (d+\beta\gamma_t) f_t-\frac{1}{2} {f_t}^2 - \frac{1}{2} (h_t - \bar{h})^2  \\
&+& - \frac{1}{2} (1-\Psi_t)^2 +  \lambda_t [\delta (1-\Psi_t)-\eta f_t\gamma_t] +\mu_t(bf_t -g_t).
\end{eqnarray*}
Let $\gamma^*_t= bang[1,0;\lambda_t- \frac{\beta}{\eta}]$ (see \eqref{eq: gamma lin}). According to the first-order conditions from the Maximum Principle (Theorem~\ref{Maximum principle}), the optimal solution  $\vec{e^*} (\gamma_t^*)= \vec{e^*}:=(h_t^*, \Psi_t^*,f_t^*, g_t^*,\gamma_t^*,  \lambda_t, \mu_t)$\footnote{For simplicity, we omit the explicit dependence of the steady state and optimal paths on $\gamma_t^*$ in the notation, except when it is necessary for the calculations.}  satisfies:
\begin{eqnarray}
\label{FOC}
 \frac{\partial H_1}{\partial f_t}(\vec{e^*})=0 &\iff& d+\beta \gamma_t^*- f_t^* -\lambda_t \eta \gamma_t^*+ \mu_t b = 0 \nonumber \\
 \frac{\partial H_1}{\partial g_t}(\vec{e^*})=0 &\iff&  h_t^* -g_t^* -\mu_t =0 \nonumber \\
\dot{\mu}_t= \rho\mu_t  -\frac{\partial H_1}{\partial h_t}(\vec{e^*}) &\iff&   \rho \mu_t - \dot{\mu}_t=   g_t^*-h_t^* + \bar{h} \nonumber \\
 \dot{\lambda_t}= \rho\lambda_t -\frac{\partial H_1}{\partial \Psi_t}(\vec{e^*})&\iff&   \rho \lambda_t - \dot{\lambda}_t = 1 - \Psi_t^* - \lambda_t \delta \\
 &&\dot{h}_t^*= bf_t^* -g_t^* \nonumber \\
 &&\dot{\Psi}_t^*=\delta (1-\Psi^*_t)-\eta \gamma_t^* f_t^* \nonumber
\end{eqnarray}

Substituting   $\delta= 1 - \rho$,  this leads to the following dynamical system


\begin{equation} \begin{pmatrix} \dot{h_t}\\ \Dot{\Psi_t}\\ \dot{\mu_t}\\ \Dot{\lambda_t} \\ \end{pmatrix} = \begin{pmatrix} b   (d+\beta\gamma_t^*) \\ 1- \rho -\eta \gamma_t^*  (d+\beta\gamma_t^*) \\ -\bar{h} \\ -1 \\ \end{pmatrix} + \begin{pmatrix} -1 & 0                  & b^2+1             & -\gamma_t^*b \eta   \\ 0  &  \rho-1        & - \eta \gamma_t^*b          & (\gamma_t^*\eta)^2    \\ 0  & 0                  & \rho+1            & 0       \\  0  & 1                  & 0                & 1         \\
\end{pmatrix}   \begin{pmatrix}
h_t^*\\ \Psi_t^* \\ \mu_t^*\\ \lambda_t^* \\ \end{pmatrix} \label{eq: system}. \end{equation}

The steady state $(h^{e}(\gamma_t^*), \Psi^{e}(\gamma_t^*),\mu^{e}(\gamma_t^*), \lambda^{e}(\gamma_t^*))=(h^{e}, \Psi^{e},\mu^{e}, \lambda^{e} )$ of \eqref{eq: system}, given by 

\begin{equation}
\begin{aligned}
h^{e}(\gamma_t^*)= &h^{e} =g^{e} +\mu^{e}= \frac{ \frac{\bar{h}}{\rho + 1}   (1+ b^2+\frac{(\eta\gamma_t^*)^2}{1-\rho})   + b   (\beta\gamma_t^* + d)}{ \frac{(\eta\gamma_t^*)^2} { 1 - \rho}+1}>0 \\
\Psi^{e}(\gamma_t^*) =&\Psi^{e} = 1- \frac{ \eta\gamma_t^* (  \frac{\bar{bh}}{\rho + 1}+\beta \gamma_t^*+ d )}{(\eta\gamma_t^*)^2 + 1 - \rho}\le 1 \\
\mu^{e} (\gamma_t^*)=&\mu^{e} = \frac{\bar{h}}{\rho + 1} \in (0, 1)\\
\lambda^{e}(\gamma_t^*) =&\lambda^{e}= 1- \Psi^{e} = \frac{ \eta\gamma_t^* (  \frac{\bar{bh}}{\rho + 1}+\beta \gamma_t^*+ d )}{(\eta\gamma_t^*)^2 + 1 - \rho}  \ge 0, \\
\text{where }
g^{e}(\gamma_t^*) =& g^{e}= b f^{e} = b\frac{ \frac{b \bar{h}}{\rho + 1} + \beta\gamma_t^* +d }{\frac{\eta^2\gamma_t^*}{1-\rho}+1} \\
\text{and  }
f^{e}(\gamma_t^*) =&f^{e}= \frac{1-\rho}{\eta}\lambda_e=\frac{ \frac{b \bar{h}}{\rho + 1} + \beta \gamma_t^* +d }{\frac{\eta^2\gamma_t^*}{1-\rho}+1},
\end{aligned}
\label{eqB: eq points}
\end{equation}
is a saddle point,  as the associated matrix has two positive and two negative eigenvalues. Thus, the path leading to the efﬁcient steady state is unstable, except along the stable trajectories whose directions are given by the eigenvectors:
$$ \vec{v_1}=
\begin{pmatrix}
    1\\
    0 \\
    0\\
    0
\end{pmatrix} , \quad \vec{v_2}= \begin{pmatrix}
    \frac{2 b \eta \gamma_t^*}{ 1+\theta_2}\\
    2 (1-\theta_2) \\
    0\\
    -2
\end{pmatrix} $$
associated to the negative eigenvalues, $\theta_1=-1$ and  

\begin{equation}
\theta_2=\theta_2(\gamma_t^*)=\frac{\rho-\sqrt{(2-\rho)^2+4(\gamma_t^*\eta)^2}}{2}=
\left\{
\begin{aligned}
 \frac{\rho-\sqrt{(2-\rho)^2+4                     \eta^2}}{2} \text{ if }& \gamma_t^*=1\\
\rho -1 \quad \quad \quad \quad  \text{ if }& \gamma_t^*=0.
 \end{aligned}
\label{theta2}
\right. 
\end{equation}
Only these stable intertemporal trajectories are admissible from an economic point of view, since they ensure that small deviations from initial conditions do not lead the system away from optimal behavior.

The stable intertemporal efﬁcient paths satisfying the initial conditions $h_0=1$ and $\Psi_0=1$ are:  
\begin{equation}
\begin{aligned}
h_t^* &= h^{e} + c_1   e^{-t} + c_2 \frac{2 b \eta \gamma_t^*}{1+\theta_2} e^{ \theta_2  t}, \\
\Psi_t^*&= \Psi^{e} +  2c_2 (1-\theta_2 ) e^{t \theta_2}, \\
\mu_t &= \mu^{e},\\
\end{aligned}
\hspace{0.7cm}
\label{eqB: optimal path}
\begin{aligned}
\lambda_t^* &= \lambda^{e} - 2 c_2    e^{ \theta_2 t}, \\
f_t^* &= d + \beta\gamma_t^* - \lambda_t \gamma_t^* \eta + \mu^{e} b, \\
g_t^* &= h_t - \mu^{e},
\end{aligned}
\end{equation}
where the constants $c_1 = 1- h^{e} -   \frac{b \eta \gamma_t^* \lambda^{e}}{1-\theta_2^2}$ and 
$c_2 =\frac{1- \Psi^{e}}{ 2(1-\theta_2)}= \frac{\lambda^{e}}{ 2(1-\theta_2)}$ are derived from  the  initial conditions. These trajectories satisfy the transversality conditions. 

To determine the optimal control $\gamma_t^*$ in accordance with the bang-bang solution framework, we analyze the conditions under which $\gamma_t^*=0$ or $\gamma_t^*=1$ satisfies the necessary optimality criteria. 

First, consider the case where  $\gamma_t^*=1$.  Using the  relations   $h^{e}(1)= \mu^{e}(1) +\frac{b\delta}{\eta}\lambda^{e}(1)$ and $\lambda^{e}(1)= \frac{\eta(b\mu^{e}+\beta+d)}{\eta^2 +\delta}$, it can be shown that the condition  \eqref{linealparamrestric} ensures both $h^{e}(1)$ and $\lambda^{e}(1)$ are less than one. Consequently, $h^{e}(1),\Psi^{e}(1) \in [0,1]$. Notice that, in the specific case where $\delta b>\eta$, the minimum in \eqref{linealparamrestric} is given by $m_1$.

If $\gamma_t^*=0$, the steady states $(h^{e}(0),f^{e}(0),g^{e}(0),\mu^{e}(0)) $ coincide with those presented in \cite{RefWorks:RefID:5-martin2013potential}. 

Now, observe that when $\gamma^*_t=1$, the adjoint operator $\lambda_t(1)$ as defined by \eqref{eqB: optimal path} is 
an increasing function of $t \in (0, \infty)$ with 
$$\lambda_0(1)= \frac{\theta_2(1)}{\theta_{2}(1) -1}\lambda^{e}(1) <\lambda^{e}(1) \  \text{ and } \ \lim_{t\to \infty}\lambda_t(1)=\lambda^e(1).$$ 

If  $\beta> \eta$, then $\lambda_t (1) \le \lambda^{e}(1) \le \frac{\beta}{\eta}$ for all $t \in (0, \infty)$, because \eqref{linealparamrestric} guarantees that $\lambda_e(1) \le 1$. In this case, the stable intertemporal efficient paths are given by \eqref{eqB: optimal path} with $\gamma_t^* = 1$.
 
Finally, if  $\beta< \eta$, it can be verified that
\begin{equation}
 \label{keyineq}
 \lambda^{e}(1)\le \frac{\beta}{\eta} \iff f^e(0)=\frac{b\bar{h}}{2- \delta}+d\le\frac{\beta \delta}{\eta^2}
\end{equation}

This concludes the proof of  Proposition~\ref{propgamma=1}.

Second-order conditions are also satisfied because the standard Hamiltonian
 
\begin{equation}
\label{eqB: current value second order}
\begin{aligned}
& \hat{H}_1( h_t, \Psi_t, \lambda_t, \mu_t) = e^{-t\rho} (-\frac{1}{2}{g_t^*}^2 + g_t^* h_t + (d+\beta {\gamma_t^*})f_t^* -\frac{1}{2} {f_t^*}^2 - \frac{1}{2} (h_t - \bar{h})^2  - \frac{1}{2} (1-\Psi_t)^2 )\\
&+ \lambda_t [\delta(1-\Psi_t) -\eta f_t^* \gamma_t^*] + \mu_t [b f_t^* - g_t^*]
\end{aligned}
\end{equation}
is concave in $( h_t, \Psi_t).$    
\end{proof}

\subsection{Proof of Proposition~\ref{prop:concave_equilibrium}}
\label{appendix D}
\setcounter{equation}{0}
\setcounter{figure}{0}
\begin{proof}
Suppose that $0<\beta\le 2\eta$  and that  the following parametric restriction holds:
\begin{equation}
   \label{paramrestr} 
   \frac{\beta (1-\rho)}{2\eta^2}-\frac{\beta}{2}\le  b\frac{\bar{h}}{\rho+1}+d \le \min \big\{\frac{4\eta (1-\rho)}{\beta^2}-\frac{\beta^2}{4\eta},\frac{1}{b}\big( 1-\frac{\bar{h}}{\rho+1}\big) -\frac{\beta}{2} \big\}
\end{equation}
with $1-\rho > \frac{\beta^4}{16\eta^2}$ and $\frac{\bar{h}}{\rho+1}\le 1-\frac{b\beta}{2}$.
We define the current-value Hamiltonian associated with the optimal control problem \eqref{eq:optimal control} with $D(\gamma_t)= \beta \gamma_t^{\frac{1}{2}}$ as:

\begin{equation*}
\begin{aligned}
& H(h_t, \Psi_t,f_t,g_t, \gamma_t, \lambda_t, \mu_t) = -\frac{1}{2}g_t^2 +g_t h_t + (d+\beta \gamma_t^{\frac{1}{2}})f_t -\frac{1}{2} f_t^2  - \frac{1}{2} (h_t - \bar{h})^2  - \frac{1}{2} (1-\Psi_t)^2  + \\ 
 & \lambda_t [\delta(1-\Psi_t) -\eta f_t \gamma_t]  + \mu_t [b f_t -g_t].
\end{aligned}
\label{eqD: current value}
\end{equation*}
The efficient solution $\vec{e^{*}}:=(h_t^{*}, \Psi_t^{*}, f_t^{*}, g_t^{*},\gamma_t^{*}, \lambda_t, \mu_t)$ satisfies the following first-order conditions of the Maximum Principle (Theorem~\ref{Maximum principle}):

\begin{eqnarray}
\label{FOC 2}
 H_{f_t}(\vec{e^{*}})=0 &\iff& d+\beta {\gamma_t^{*}}^{\frac{1}{2}}- f_t^{*} -\lambda_t \eta \gamma_t^{*}+ \mu_t b = 0 \nonumber \\
 H_{g_t}(\vec{e^{*}})=0 &\iff&  h_t^{*} -g_t^{*} -\mu_t =0 \nonumber \\
  H_{\gamma_t}(\vec{e^{*}})=0 &\iff& \frac{1}{2} \beta  f_t^{*} {\gamma_t^{*}}^{-\frac{1}{2}} - \lambda_t \eta f_t^{*}=0  \nonumber \\
\dot{\mu}_t= \rho\mu_t  -H_{h_t}(\vec{e^{*}}) &\iff&   \rho \mu_t - \dot{\mu}_t=   g_t^{*}-h_t^{*} + \bar{h} \nonumber \\
 \dot{\lambda_t}= \rho\lambda_t -H_{\Psi_t}(\vec{e^{*}})&\iff&   \rho \lambda_t - \dot{\lambda}_t = 1 - \Psi_t^{*} - \lambda_t \delta \\
 &&\dot{h}_t^{*}= bf_t^{*} -g_t^{*} \nonumber \\
 &&\dot{\Psi}_t^{*}=\delta (1-\Psi^{*}_t)-\eta \gamma_t^{*} f_t^{*}- \nonumber
\end{eqnarray}

By setting  $\delta= 1 - \rho$, we obtain  the following system of equations:
\begin{eqnarray}
\label{eqB: system two}
\begin{array}{cc}
g_t^* = h_t^* - \mu_t  & \Dot{h}_t = b  \left(d + \frac{\beta^2}{4 \eta \lambda_t^*} + b \mu_t\right) - h_t^* + \mu_t \\
\gamma_t^*= \frac{\beta^2}{4 \eta^2 \lambda_t^2}  & \Dot{\Psi}_t^* = (1-\rho)   (1- \Psi_t^*) - \eta   \left(d + \frac{\beta^2}{4 \eta \lambda_t} + b \mu_t \right) \left(\frac{\beta^2}{4 \eta^2 \lambda_t^2} \right)  \\
f_t^* = d + \frac{\beta^2}{4 \lambda_t \eta} + b \mu_t  & \Dot{\mu}_t = (1+ \rho)   \mu_t - \Bar{h} \\
& \Dot{\lambda}_t = \lambda_t + \Psi_t^*-1
\end{array}
\end{eqnarray}

The steady state  solution $(h_e, \Psi_e,\mu_e, \lambda_e) $ of \eqref{eqB: system two} is  given by: 

\begin{eqnarray}
\begin{array}{cc}
\mu_e = \frac{ \Bar{h}}{\rho + 1} &  h_e= bd + \frac{b \beta}{2}{\gamma_e}^{1/2}+ (1+b^2) \mu_e \\
 \gamma_e =\frac{\beta^2}{4 \eta^2 {\lambda_e}^2} \iff \lambda_e = \frac{\beta}{2 \eta}  {\gamma_e}^{-1/2} &  g_e = h_e - \mu_e \\
\psi_e = 1 - \lambda_e  &  f_e = d + \frac{\beta^2}{4 \lambda_e \eta} + b \mu_e\\
  (d+ \frac{\beta}{2}{\gamma_e}^{1/2} + \mu_e b ) \frac{2 \eta^2}{\beta (1- \rho)} {\gamma_e}^{3/2} -1=0 .& \\
\label{eqD: eq points}
\end{array}
\end{eqnarray}

Given that $\Psi_e =1-\lambda_e$ and $\lambda_e=\frac{\beta}{2\eta}\gamma_e^{-1/2}$, the conditions $\Psi_e, \lambda_e \in [0,1]$ hold if  $\frac{\beta^2}{4\eta^2} \le \gamma_e \le 1$ with $\beta>0$.

Now, define the function 
\begin{equation}\label{eqD: f. gamma}
    P(\gamma_e) :=(d+ \frac{\beta}{2}{\gamma_e}^{1/2} + \mu_e b ) \frac{2 \eta^2}{\beta (1- \rho)} {\gamma_e}^{3/2} -1.
\end{equation}
On the one had, the parametric restrictions \eqref{paramrestr} ensure that  $P\left(\frac{\beta^2}{4\eta^2}\right) \le 0 \le P(1)$. On the other hand,  $P(\gamma_e)$ is continuous and increasing on $[0,1]$, hence  the Intermediate Value Theorem  implies that $\gamma_e \in \left[\frac{\beta^2}{4\eta^2}, 1\right]$. Finally, it can be checked that under parametric restrictions given by \eqref{paramrestr}  the restriction $h_e \le 1$ holds true.

Stability around the steady state is analyzed by linearizing the system and computing the Jacobian matrix:
\begin{equation*}
J(h_e, \Psi_e, \mu_e, \lambda_e) 
=
\begin{pmatrix}
-1 & 0 & 1+b^2 & - b \eta \gamma_e \\
0 & -(1-\rho) & - b \eta \gamma_e & (\eta\gamma_e)^2+2(1-\rho) \\
0 & 0 & 1 + \rho & 0 \\
0 & 1 & 0 & 1 \\
\end{pmatrix}.
\label{eqD: Jacb}
\end{equation*}

The eigenvalue structure (two positive and two negative roots) confirms that the steady state is a saddle point. This permits the construction of a two-dimensional stable manifold, along which the system converges to equilibrium from given initial conditions.

The Jacobian matrix provides crucial information about the local stability of an equilibrium point. By computing its  eigenvalues, one can classify the nature of the steady state. In this case, the matrix has two positive and two negative eigenvalues, indicating that the equilibrium is a saddle point.
According to Hartman’s linearization theorem \citep{hartmant}, if all eigenvalues of the Jacobian matrix have non-zero real parts—i.e., the equilibrium is hyperbolic—then the qualitative behavior of the nonlinear system in a neighborhood of that point is locally equivalent to that of the linearized system. Therefore, the system inherits the saddle structure of the linearized dynamics, and a two-dimensional stable manifold exists. This manifold corresponds to the span of the eigenvectors associated with the negative eigenvalues, and any initial condition lying within it leads the system to converge to the equilibrium. Outside this manifold, trajectories diverge.

The eigenvectors corresponding to the negative eigenvalues, $\theta_1=-1$ and $\theta_2= \frac{\rho - \sqrt{(\rho-2)^2+4W}}{2}$, are given by:
$$ \vec{v_1}=
\begin{pmatrix}
    1\\
    0 \\
    0\\
    0
\end{pmatrix} , \quad \vec{v_2}= \begin{pmatrix}
    \frac{-b\eta \gamma_e}{1+\theta_2}\\
    \frac{W}{1-\rho+\theta_2} \\
    0   \\
    1
\end{pmatrix}, $$
where $W=(\eta \gamma_e)^2+2(1-\rho)$.

Therefore,  the efficient intertemporal paths  satisfying the initial conditions $h_0=1$ and $\Psi_0=1$ are given by:

\begin{equation}
\left\{
\begin{aligned}
h_t^{*} &= h_e + c_1   e^{-t} - \frac{ b \eta  \gamma_e}{1+\theta_2} c_2 e^{t \theta_2}\\
\Psi_t^{*}&= \Psi_e + \frac{ W}{1-\rho+\theta_2} c_2 e^{t\theta_2}   \\
\mu_t &= \mu_e \\
\lambda_t &= \lambda^{e}(1) + c_2 e^{t \theta_2}\\
f_t^{*} &= d + \frac{\beta^2}{4 \lambda_t \eta} + b \mu_e\\
g_t^{*} &= h_t - \mu_e \\
\gamma_t^{*} &= \frac{\beta^2}{4 \eta^2 \lambda_t^2} \\
\end{aligned}
\label{eqD: optimal path}
\right. 
\implies 
\left\{ \begin{aligned}
h_t^{*} &= h_e+  (1  - h_e) e^{-t}- \frac{ b\eta \gamma_e}{1+\theta_2} c_2 \big(e^{ \theta_2 t}-e^{-t}\big) \\
\Psi_t^{*}&= \Psi_e +  (1 - \Psi_e) e^{t \theta_2} \\
\mu_t &= \mu_e \\
\lambda_t &= \lambda_e  \big(1+ \frac{1-\rho+\theta_2}{W} e^{ \theta_2 t}\big) \\
f_t^{*} &= d + \frac{\beta^2}{4 \lambda_t \eta} + b \mu_e \\
g_t^{*} &= h_t - \mu_e \\
\gamma_t^{*} &= \frac{\beta^2}{4 \eta^2 \lambda_t^2} \\
\end{aligned}
\right.
\end{equation}

where $c_1 =1 - h_e + c_2\frac{ b \eta \gamma_e}{1+ \theta_2}$ and 
$c_2 = \frac{(1-\rho+\theta_2)(1-\Psi_e)}{W}$ are obtained from  the  initial conditions $h_0=1$ and $\Psi_0=1$.

Regarding the control variable $\gamma_t^* = \frac{\beta^2}{4 \eta^2 \lambda_t^2}$, which must lie within $[0,1]$, we enforce this constraint explicitly using a Karush-Kuhn-Tucker (KKT) formulation:
\begin{equation*}
\gamma_t^* = \min\left\{1,\,\max\left\{0, \frac{\beta^2}{4\eta^2\lambda_t^2} \right\}\right\}.
\end{equation*}

The transversality conditions $\lim_{t \to \infty} e^{-\rho t} \mu_t=0 $ and $\lim_{t \to \infty} e^{-\rho t} \lambda_t=0$ are satisfied, and the second-order conditions hold as the  Standart Hamiltonian:
 
\begin{equation}
\label{eqD: current value second order}
\begin{aligned}
& \hat{H}( h_t, \Psi_t, \lambda_t, \mu_t) = e^{-t\rho} (-\frac{1}{2}{g_t^{*}}^2 + g_t^{*} h_t + (d+\beta {\gamma_t^{*}}^{1/2})f_t^{*} -\frac{1}{2} {f_t^{*}}^2 - \frac{1}{2} (h_t - \bar{h})^2  - \frac{1}{2} (1-\Psi_t)^2 )\\
&+ \lambda_t [\delta(1-\Psi_t) -\eta f_t^{*} \gamma_t^{**}] + \mu_t [b f_t^{*} - g_t^{*}]
\end{aligned}
\end{equation}
is concave in $( h_t, \Psi_t).$

Finally, as $\beta \to 0$, the optimal fertilizer use converges to zero: 
\begin{equation}
\lim_{\beta \to 0} \frac{\beta^2}{4  \eta^2  \lambda^2} = 0.
\label{eqD: demo}
\end{equation}
Accordingly, the optimal solution coincides with the benchmark model without fertilizers proposed by \cite{RefWorks:RefID:5-martin2013potential}.
\end{proof}

\end{document}